\documentclass{svmult}
\usepackage{bbm}
\usepackage{amsfonts}
\usepackage{mathrsfs}
\usepackage{amssymb}
\newenvironment{enumerate-roman}{\begin{enumerate}}{\end{enumerate}}
\usepackage{makeidx}         
\usepackage{graphicx}        
\usepackage{multicol}        
\usepackage[bottom]{footmisc}

\makeindex             

\newtheorem{thm}{Theorem}[section]
\newtheorem{lem}[thm]{Lemma}
\newtheorem{prop}[thm]{Proposition}
\newtheorem{cor}[thm]{Corollary}
\newtheorem{rmk}[thm]{Remark}

\newtheorem{defi}[thm]{Definition}

\begin{document}

\baselineskip13pt

\title*{Stationary Solutions of SPDEs and Infinite Horizon BDSDEs with Non-Lipschitz Coefficients}
\titlerunning{Stationary Solution of SPDEs with Non-Lipschitz Coefficients}

\author{Qi Zhang\inst{\ {\rm a}, {\rm b}, {\rm c}}, Huaizhong Zhao\inst{\ {\rm a}}}
\authorrunning{Q. Zhang and H.Z. Zhao}
\institute{$^{\rm a}$ Department of Mathematical Sciences,
Loughborough University, Loughborough, LE11 3TU, UK.\\
$^{\rm b}$ School of Mathematics and System Sciences, Shandong
University, Jinan, 250100, China.\\
$^{\rm c}$ Current address: School of Mathematical Sciences, Fudan
University, Shanghai, 200433, China.\\\texttt{Emails:
qzh@fudan.edu.cn}; \texttt{H.Zhao@lboro.ac.uk}}

\maketitle
\newcounter{bean}
\begin{abstract}
We prove a general theorem that the
$L_{\rho}^2({\mathbb{R}^{d}};{\mathbb{R}^{1}})\otimes
L_{\rho}^2({\mathbb{R}^{d}};{\mathbb{R}^{d}})$ valued solution of an
infinite horizon backward doubly stochastic differential equation,
if exists, gives the stationary solution of the corresponding
stochastic partial differential equation. We prove the existence and
uniqueness of the
$L_{\rho}^2({\mathbb{R}^{d}};{\mathbb{R}^{1}})\otimes
L_{\rho}^2({\mathbb{R}^{d}};{\mathbb{R}^{d}})$ valued solutions for
backward doubly stochastic differential equations on finite and
infinite horizon with linear growth without assuming Lipschitz
conditions, but under the monotonicity condition. Therefore the
solution of finite horizon problem gives the solution of the initial
value problem of the corresponding stochastic partial differential
equations, and the solution of the infinite horizon problem gives
the stationary solution of the SPDEs according to our general
result.
\end{abstract}
\vskip5pt

\vskip5pt

\renewcommand{\theequation}{\arabic{section}.\arabic{equation}}

\section{Introduction, Basic Notation and Main Results}\label{s31}
The notion of the pathwise stationary solutions for stochastic
partial differential equations (SPDEs) is a fundamental concept in
the study of the long time behaviour of the stochastic dynamical
systems driven by the SPDEs. It describes the pathwise invariance of
the stationary solution, over time, along the measurable and
$P$-preserving transformation $\theta_t:\Omega\to\Omega$ and the
pathwise limit of the solutions of the random dynamical systems. It
is a random fixed point $Y(\omega)$ in the state space of the random
dynamical system, in the sense that the solution
$v(t,Y(\omega),\omega)$ of the SPDE with the initial value
$Y(\omega)$ is  equal to $Y(\theta_t\omega)$, which is still $Y$,
but with a different sample path $\theta_t\omega$. Therefore
$Y(\theta_t\omega)$ is a particular solution of the SPDE with the
pathwise stationary property. Needless to say that the ``one-force,
one-solution" setting is a natural extension of the equilibrium or
fixed point in the theory of the deterministic dynamical systems to
stochastic counterparts. Such a random fixed point consists of
infinitely many randomly moving invariant surfaces on the
configuration space due to the random external force pumped to the
system constantly. The study of its existence and stability is of
great interests in both mathematics and physics. We would like to
point out that the existence of stationary solutions is a basic
assumption in many works on random dynamical systems e.g. in the
study of stability (Has$'$minskii \cite{ha}), and in the theory of
stable and unstable manifolds (Arnold \cite{ar}, Mohammed, Zhang and
Zhao \cite{mo-zh-zh}, Duan, Lu and Schmalfuss \cite{du-lu-sc1}). But
these theories give neither the existence of stationary solutions,
nor a way of finding them. However, in contrast to the deterministic
dynamical systems, the existence of stationary solutions of random
dynamical systems is a more difficult and subtle problem. It is easy
to see that the solutions of elliptic type partial differential
equations give the stationary solutions of the corresponding
parabolic type partial differential equations, though the elliptic
partial differential equations are difficult problems to study as
well. However, for stochastic partial differential equations of the
parabolic type, such kind of connection does not exist. In
\cite{mo-zh-zh}, Mohammed, Zhang and Zhao introduced an integral
equation of infinite horizon for the stationary solutions of certain
stochastic evolution equations. But the existence of the solutions
of such stochastic integral equations in general is far from clear.
In \cite{zh-zh}, Zhang and Zhao proved that the solution of an
infinite horizon backward doubly stochastic differential equation
(BDSDE), if exists, is a perfect stationary solution. Moreover,
under the Lipschitz and monotone conditions, the solution indeed
exists and gives the stationary solution of the corresponding SPDEs
of the parabolic type. It was known that the solutions of infinite
horizon backward stochastic differential equations (BSDEs) give a
classical or viscosity solution of elliptic type partial
differential equations (Poisson equations) from the works of Peng
\cite{pe} and Pardoux \cite{pa}. So philosophically, it is very
natural to represent the stationary solutions of SPDEs as the
solutions of the infinite horizon BDSDEs, like the case of the
Poisson equations as the solutions of the infinite horizon backward
stochastic differential equations. Other works on stationary
solutions of certain types of SPDEs usually under additive or linear
noise include Sinai \cite{si1}, \cite{si2}, Caraballo, Kloeden,
Schmalfuss \cite{kloeden}.

In this paper, we will put above idea on infinite horizon BDSDEs in
a general setting and prove a general theorem which basically says,
if the infinite horizon BDSDE has a unique solution in the space
$S^{2,-K}\bigcap
M^{2,-K}([0,\infty);L_{\rho}^2\\({\mathbb{R}^{d}};{\mathbb{R}^{1}}))\bigotimes
M^{2,-K}([0,\infty);L_{\rho}^2({\mathbb{R}^{d}};{\mathbb{R}^{d}}))$
for a $K>0$, and the finite horizon BDSDE gives the representation
for the solution of the corresponding SPDE, then, the solution of
the infinite horizon BDSDE gives the stationary solution of the
corresponding SPDE. Following this result, to study the existence of
stationary solutions of SPDEs is transformed to study the existence
and uniqueness of the solutions of the corresponding infinite
horizon BDSDEs. In \cite{zh-zh}, we studied such equations when the
nonlinear coefficients are assumed to be Lipschitz continuous. In
this paper, we continue our work \cite{zh-zh} to study the weak
solution (in the weighted Sobolev space
$H^{1}_{\rho}({\mathbb{R}^{d}};{\mathbb{R}^{1}})$ space) of the
following parabolic SPDE without assuming the Lipschitz continuity
of $f$ on $v$:
\begin{eqnarray}\label{zhao21}
dv(t,x)&=&[\mathscr{L}v(t,x)+f\big(x,v(t,x),\sigma^*(x)Dv(t,x)\big)]dt\nonumber\\
&&+g\big(x,v(t,x), \sigma^*(x)Dv(t,x)\big)d{B}_t.
\end{eqnarray}
Here ${B}$ is a two-sided cylindrical Brownian motion valued on a
separable Hilbert space $U_0$ in a probability space $(\Omega,
\mathscr{F}, P)$; $\mathscr{L}$ is the infinitesimal generator of a
diffusion process $X_{s}^{t,x}$ (the solution of Eq.(\ref{qi17}))
given by
\begin{eqnarray}\label{buchong1}
\mathscr{L}={1\over2}\sum_{i,j=1}^na_{ij}(x){{\partial^2}\over{\partial
x_i\partial x_j}}+\sum_{i=1}^nb_i(x){\partial\over{\partial x_i}}
\end{eqnarray}
with $\big(a_{ij}(x)\big)=\sigma\sigma^*(x)$;
$L_{\rho}^2({\mathbb{R}^{d}};{\mathbb{R}^{1}})$ is the Hilbert space
with the inner product
\begin{eqnarray*}
\langle
u_1,u_2\rangle=\int_{\mathbb{R}^d}u_1(x)u_2(x)\rho^{-1}(x)dx,
\end{eqnarray*}
i.e. a $\rho$-weighted $L^2$ space, where $\rho(x)=(1+|x|)^q$,
$q>3$, is a weight function. It is easy to see that
$\rho(x):\mathbb{R}^d\longrightarrow\mathbb{R}^1$ is a continuous
positive function satisfying
$\int_{\mathbb{R}^{d}}|x|^p\rho^{-1}(x)dx<\infty$ for any
$p\in(2,q-1)$. Note that we can consider more general $\rho$ which
satisfies the above condition and conditions in \cite{ba-ma} and all
the results of this paper still hold. The SPDEs we study in this
paper are very general with the noise term $g$ being allowed to be
nonlinear in $v$ and $\nabla v$. However, many techniques of
\cite{zh-zh} when $f$ is Lipschitz do not work here. Although we can
follow a similar procedure, as in \cite{zh-zh}, to consider first
the finite horizon BDSDEs, then to make the connection with the weak
solutions of the corresponding SPDEs and to find a Cauchy sequence
to pass the terminal time of BDSDEs to infinity, we have to
introduce new techniques to deal with the difficulties arising from
the lack of the Lipschitz continuity of $f$ on $y$.

Define $u(t,x)=v(T-t,x)$ for arbitrary $T$ and $0\leq t\leq T$. We
can show that $u$ satisfies the following backward SPDE:
\begin{eqnarray}
\left\{\begin{array}{l}\label{zhang685}
du(t,x)+[\mathscr{L}u(t,x)+f\big(x,u(t,x),(\sigma^*\nabla
u)(t,x)\big)]dt\\
\ \ \ \ \ \ \ \ \ \ -g\big(x,u(t,x),(\sigma^*\nabla
u)(t,x)\big)d^\dagger
\hat{B}_t=0\\
u(T,x)=v(0,x).
\end{array}\right.
\end{eqnarray}
Here $\mathscr{L}$ is given by (\ref{buchong1}) and
$\hat{B}_s=B_{T-s}-B_T$. Let ${\cal N}$ denote the class of $P$-null
sets of ${\mathscr{F}}$. We define
\begin{eqnarray*}
&&\mathscr{F}_{t,T}\triangleq{\mathscr{F}_{t,T}^{\hat{B}}}\bigvee
\mathscr{F}_t^W\bigvee{\cal N},\ \ \ {\rm for}\ 0\leq t\leq T;\\
&&\mathscr{F}_t\triangleq{\mathscr{F}_{t,\infty}^{\hat{B}}}\bigvee
\mathscr{F}_t^W\bigvee{\cal N},\ \ \ \ \ {\rm for}\ t\geq0.
\end{eqnarray*}
Recall Definitions 2.1 and 2.2 in \cite{zh-zh}:
\begin{defi}\label{qi00}
Let $\mathbb{S}$ be a Hilbert space with norm $\|\cdot\|_\mathbb{S}$
and Borel $\sigma$-field $\mathscr{S}$. For $K\in\mathbb{R}^+$, we
denote by $M^{2,-K}([0,\infty);\mathbb{S})$ the set of
$\mathscr{B}_{\mathbb{R}^+}\otimes\mathscr{F}/\mathscr{S}$
measurable random processes $\{\phi(s)\}_{s\geq0}$ with values in
$\mathbb{S}$ satisfying
\begin{enumerate-roman}
\item $\phi(s):\Omega\rightarrow\mathbb{S}$ is $\mathscr{F}_s$ measurable for $s\geq 0$;
\item $E[\int_{0}^{\infty}{\rm e}^{-Ks}\|\phi(s)\|_\mathbb{S}^2ds]<\infty$.
\end{enumerate-roman}
Also we denote by $S^{2,-K}([0,\infty);\mathbb{S})$ the set of
$\mathscr{B}_{\mathbb{R}^+}\otimes\mathscr{F}/\mathscr{S}$
measurable random processes $\{\psi(s)\}_{s\geq0}$ with values in
$\mathbb{S}$ satisfying
\begin{enumerate-roman}
\item $\psi(s):\Omega\rightarrow\mathbb{S}$ is $\mathscr{F}_s$ measurable for $s\geq0$ and $\psi(\cdot,\omega)$ is continuous $P$-a.s.;
\item $E[\sup_{s\geq0}{\rm
e}^{-Ks}\|\psi(s)\|_\mathbb{S}^2]<\infty$.
\end{enumerate-roman}
\end{defi}
\begin{defi}\label{zhao005}
Let $\mathbb{S}$ be a Hilbert space with norm $\|\cdot\|_\mathbb{S}$
and Borel $\sigma$-field $\mathscr{S}$. We denote by
$M^{2,0}([t,T];\mathbb{S})$ the set of
$\mathscr{B}_{[t,T]}\otimes\mathscr{F}/\mathscr{S}$ measurable
random processes $\{\phi(s)\}_{t\leq s\leq T}$ with values in
$\mathbb{S}$ satisfying
\begin{enumerate-roman}
\item $\phi(s):\Omega\rightarrow\mathbb{S}$ is $\mathscr{F}_{s,T}\bigvee{\mathscr{F}_{T,\infty}^{\hat{B}}}$ measurable for $t\leq s\leq T$;
\item $E[\int_{t}^{T}\|\phi(s)\|_\mathbb{S}^2ds]<\infty$.
\end{enumerate-roman}
Also we denote by $S^{2,0}([t,T];\mathbb{S})$ the set of
$\mathscr{B}_{[t,T]}\otimes\mathscr{F}/\mathscr{S}$ measurable
random processes $\{\psi(s)\}_{t\leq s\leq T}$ with values in
$\mathbb{S}$ satisfying
\begin{enumerate-roman}
\item $\psi(s):\Omega\rightarrow\mathbb{S}$ is
$\mathscr{F}_{s,T}\bigvee{\mathscr{F}_{T,\infty}^{\hat{B}}}$
measurable for $t\leq s\leq T$ and $\psi(\cdot,\omega)$ is
continuous $P$-a.s.;
\item $E[\sup_{t\leq s\leq T}\|\psi(s)\|_\mathbb{S}^2]<\infty$.
\end{enumerate-roman}
\end{defi}

Recall also the weak solution of the SPDE (\ref{zhang685}) as
follows for the convenience of the reader.
\begin{defi}\label{qi042}
A process $u$ is called a weak solution {\rm(}solution in
$L_{\rho}^2({\mathbb{R}^{d}};{\mathbb{R}^{1}})${\rm)} of
Eq.(\ref{zhang685}) if $(u,\sigma^*\nabla u)\in
M^{2,0}([0,T];L_{\rho}^2({\mathbb{R}^{d}};{\mathbb{R}^{1}}))\bigotimes
M^{2,0}([0,T];L_{\rho}^2({\mathbb{R}^{d}};{\mathbb{R}^{d}}))$ and
for an arbitrary $\varphi\in
C_c^{\infty}(\mathbb{R}^d;\mathbb{R}^1)$,
\begin{eqnarray}\label{qi16}
&&\int_{\mathbb{R}^{d}}u(t,x)\varphi(x)dx-\int_{\mathbb{R}^{d}}u(T,x)\varphi(x)dx\nonumber\\
&&-{1\over2}\int_{t}^{T}\int_{\mathbb{R}^{d}}(\sigma^*\nabla u)(s,x)(\sigma^*\nabla\varphi)(x)dxds\nonumber\\
&&-\int_{t}^{T}\int_{\mathbb{R}^{d}}u(s,x)div\big((b-\tilde{A})\varphi\big)(x)dxds\nonumber\\
&=&\int_{t}^{T}\int_{\mathbb{R}^{d}}f\big(x,u(s,x),(\sigma^*\nabla u)(s,x)\big)\varphi(x)dxds\\
&&-\sum_{j=1}^{\infty}\int_{t}^{T}\int_{\mathbb{R}^{d}}g_j\big(x,u(s,x),(\sigma^*\nabla
u)(s,x)\big)\varphi(x)dxd^\dagger{\hat{\beta}}_j(s)\ \ \ P-{\rm
a.s.}\nonumber
\end{eqnarray}
Here $\tilde{A}_j\triangleq{1\over2}\sum_{i=1}^d{\partial
a_{ij}(x)\over\partial x_i}$, and
$\tilde{A}=(\tilde{A}_1,\tilde{A}_2,\cdot\cdot\cdot,\tilde{A}_d)^*$.
\end{defi}
\begin{rmk}
The weak solution of the forward SPDE (\ref{zhao21}) can be defined
similarly. Sometimes in this paper, we denote it by $v(t,v_0)(\cdot)$
to emphasize its dependence on its initial
value $v_0$.
\end{rmk}
For $k\geq0$, we denote by $C_{l,b}^k$
\ the set of $C^k$-functions whose partial derivatives of order less
than or equal to $k$ are bounded and by
$H^k_\rho$
\ the $\rho$-weighted Sobolev space (See e.g. \cite{ba-ma}). In
order to connect BDSDEs with SPDEs, the form of BDSDEs should be a
kind of FBDSDEs (forward and backward doubly SDEs). So we first
 let $X_{s}^{t,x}$ be a diffusion process given by the
solution of the following forward SDE for $s\geq t$,
\begin{eqnarray}\label{qi17}
X_{s}^{t,x}=x+\int_{t}^{s}b(X_{u}^{t,x})du+\int_{t}^{s}\sigma(X_{u}^{t,x})dW_u,
\end{eqnarray}
where $b\in C_{l,b}^2(\mathbb{R}^{d};\mathbb{R}^{d})$, $\sigma\in
C_{l,b}^3(\mathbb{R}^{d};\mathbb{R}^{d}\times\mathbb{R}^{d})$, and
for $0\leq s<t$, we regulate $X_{s}^{t,x}=x$. We now construct the
measurable metric dynamical system through defining a measurable and
measure-preserving shift. Let
$\hat{\theta}_t:\Omega\longrightarrow\Omega$, $t\geq0$, be a
measurable mapping on $(\Omega, {\mathscr{F}}, P)$, defined by
\begin{eqnarray*}
\hat{\theta}_{t}\circ \hat{B}_s=\hat{B}_{s+t}-\hat{B}_t,\ \ \
\hat{\theta}_{t}\circ W_s=W_{s+t}-W_t.
\end{eqnarray*}
Then for any $s$, $t\geq0$,
\begin{description}
\item[$(\textrm{i})$]$P\cdot\hat{\theta}_{t}^{-1}=P$;
\item[$(\textrm{i}\textrm{i})$]$\hat{\theta}_{0}=I$, where $I$ is the identity transformation on $\Omega$;
\item[$(\textrm{i}\textrm{i}\textrm{i})$]$\hat{\theta}_{s}\circ\hat{\theta}_{t}=\hat{\theta}_{s+t}$.
\end{description}
For any $r\geq0$, $s\geq t$, $x\in\mathbb{R}^d$, apply
$\hat{\theta}_r$ to SDE (\ref{qi17}), then we have
\begin{eqnarray}\label{qi18}
\hat{\theta}_r\circ X_{s}^{t,x}=X_{s+r}^{t+r,x}\ \ {\rm for}\ {\rm
all}\ r,s,t,x\ \ \ {\rm a.s.}
\end{eqnarray}

The following lemma in \cite{zh-zh} is an extension of the
equivalence of norm principle given in \cite {ku1}, \cite{ba-le},
\cite{ba-ma} to the cases when $\varphi$ and $\Psi$ are random.
\begin{lem}\label{qi045}(generalized equivalence
of norm principle \cite{zh-zh}) Let $\rho$ be the weight function
and $X$ be the diffusion process given above. If $s\in[t,T]$,
$\varphi:\Omega\times\mathbb{R}^d\rightarrow\mathbb{R}^1$ is
independent of $\mathscr{F}^W_{t,s}$ and $\varphi\rho^{-1}\in
L^1(\Omega\otimes\mathbb{R}^{d})$, then there exist two constants
$c>0$ and $C>0$ such that
\begin{eqnarray*}
cE[\int_{\mathbb{R}^{d}}|\varphi(x)|\rho^{-1}(x)dx]\leq
E[\int_{\mathbb{R}^{d}}|\varphi(X_{s}^{t,x})|\rho^{-1}(x)dx]\leq
CE[\int_{\mathbb{R}^{d}}|\varphi(x)|\rho^{-1}(x)dx].
\end{eqnarray*}
Moreover if
$\Psi:\Omega\times[t,T]\times\mathbb{R}^d\rightarrow\mathbb{R}^1$,
$\Psi(s,\cdot)$ is independent of $\mathscr{F}^W_{t,s}$ and
$\Psi\rho^{-1}\in L^1(\Omega\otimes[t,T]\otimes\mathbb{R}^{d})$,
then
\begin{eqnarray*}
&&cE[\int_{t}^{T}\int_{\mathbb{R}^{d}}|\Psi(s,x)|\rho^{-1}(x)dxds]\leq E[\int_{t}^{T}\int_{\mathbb{R}^{d}}|\Psi(s,X_{s}^{t,x})|\rho^{-1}(x)dxds]\\
&\leq&CE[\int_{t}^{T}\int_{\mathbb{R}^{d}}|\Psi(s,x)|\rho^{-1}(x)dxds].
\end{eqnarray*}
\end{lem}

We consider the following infinite horizon BDSDE:
\begin{eqnarray}\label{qi13}
{\rm e}^{-Ks}Y_{s}^{t,x}&=&\int_{s}^{\infty}{\rm e}^{-Kr}f(X_{r}^{t,x},Y_{r}^{t,x},Z_{r}^{t,x})dr+\int_{s}^{\infty}K{\rm e}^{-Kr}Y_{r}^{t,x}dr\\
&&-\int_{s}^{\infty}{\rm e}^{-Kr}
g(X_{r}^{t,x},Y_{r}^{t,x},Z_{r}^{t,x})d^\dagger{\hat{B}}_r-\int_{s}^{\infty}{\rm
e}^{-Kr}\langle Z_{r}^{t,x},dW_r\rangle.\nonumber
\end{eqnarray}
Here
$\hat{B}_r=\sum_{j=1}^{\infty}\sqrt{\lambda_j}\hat{\beta}_j(r)e_j$,
$\{\hat{\beta}_j(r)\}_{j=1,2,\cdots}$ are mutually independent
one-dimensional Brownian motions;
$f:\mathbb{R}^{d}\times\mathbb{R}^1\times\mathbb{R}^{d}{\longrightarrow{\mathbb{R}^1}}$;
$g:\mathbb{R}^{d}\times\mathbb{R}^1\times\mathbb{R}^{d}\longrightarrow
{\mathcal{L}^2_{U_0}(\mathbb{R}^{1})}$. Set $g_j\triangleq
g\sqrt{\lambda_j}e_j:\mathbb{R}^{d}\times\mathbb{R}^1\times\mathbb{R}^{d}{\longrightarrow{\mathbb{R}^1}}$,
then Eq.(\ref{qi13}) is equivalent to
\begin{eqnarray*}\label{qi14}
{\rm e}^{-Ks}Y_{s}^{t,x}&=&\int_{s}^{\infty}{\rm e}^{-Kr}f(X_{r}^{t,x},Y_{r}^{t,x},Z_{r}^{t,x})dr+\int_{s}^{\infty}K{\rm e}^{-Kr}Y_{r}^{t,x}dr\nonumber\\
&&-\sum_{j=1}^{\infty}\int_{s}^{\infty}{\rm
e}^{-Kr}g_j(X_{r}^{t,x},Y_{r}^{t,x},Z_{r}^{t,x})d^\dagger{\hat{\beta}}_j(r)-\int_{s}^{\infty}{\rm
e}^{-Kr}\langle Z_{r}^{t,x},dW_r\rangle.
\end{eqnarray*}
\begin{defi}\label{qi041} (Definition 2.7 in \cite{zh-zh})
A pair of processes $(Y_{\cdot}^{t,\cdot},Z_{\cdot}^{t,\cdot})\in
S^{2,-K}\bigcap
M^{2,-K}([0,\infty);L_{\rho}^2({\mathbb{R}^{d}};{\mathbb{R}^{1}}))\bigotimes
M^{2,-K}([0,\infty);L_{\rho}^2({\mathbb{R}^{d}};{\mathbb{R}^{d}}))$
is called a solution of Eq.(\ref{qi13}) if for an arbitrary
$\varphi\in C_c^{0}(\mathbb{R}^d;\mathbb{R}^1)$,
\begin{eqnarray}\label{qi15}
\int_{\mathbb{R}^{d}}{\rm
e}^{-Ks}Y_s^{t,x}\varphi(x)dx&=&\int_{s}^{\infty}\int_{\mathbb{R}^{d}}{\rm e}^{-Kr}f(X_{r}^{t,x},Y_r^{t,x},Z_r^{t,x})\varphi(x)dxdr\nonumber\\
&&+\int_{s}^{\infty}\int_{\mathbb{R}^{d}}K{\rm e}^{-Kr}Y_{r}^{t,x}\varphi(x)dxdr\nonumber\\
&&-\sum_{j=1}^{\infty}\int_{s}^{\infty}\int_{\mathbb{R}^{d}}{\rm e}^{-Kr}g_j(X_{r}^{t,x},Y_r^{t,x},Z_r^{t,x})\varphi(x)dxd^\dagger{\hat{\beta}}_j(r)\nonumber\\
&&-\int_{s}^{\infty}\langle \int_{\mathbb{R}^{d}}{\rm
e}^{-Kr}Z_r^{t,x}\varphi(x)dx,dW_r\rangle\ \ \ P-{\rm a.s.}
\end{eqnarray}
\end{defi}

We will prove the following theorem under a general setting.

\begin{thm}\label{qi888}
If Eq.(\ref{qi13}) has a unique solution
$(Y_{\cdot}^{t,\cdot},Z_{\cdot}^{t,\cdot})\in S^{2,-K}\bigcap
M^{2,-K}\\([0,\infty);L_{\rho}^2({\mathbb{R}^{d}};{\mathbb{R}^{1}}))
\bigotimes
M^{2,-K}([0,\infty);L_{\rho}^2({\mathbb{R}^{d}};{\mathbb{R}^{d}}))$
and $u(t,\cdot)\triangleq Y_t^{t,\cdot}$ is a continuous weak
solution of Eq.(\ref{zhang685}), then $u(t,\cdot)$ has an
indistinguishable version which is a ``perfect" stationary weak
solution of Eq.(\ref{zhang685}). Furthermore, let
$\hat{B}_s={B}_{T'-s}-{B}_{T'}$ for all $s\geq0$ in Eq. (\ref{qi13})
and $v_t(\cdot)\triangleq u(T'-t,\cdot)=Y_{T'-t}^{T'-t,\cdot}$ for
arbitrary $T'$ and $t\in[0,T']$,  then $v_t(\cdot)$ is independent
of $T'$ and is a ``perfect" stationary weak solution of
Eq.(\ref{zhao21}) i.e.
\begin{eqnarray*}
v_t(\omega)=v_0(\theta _t\omega)=v(t,v_0(\omega),\omega)\ {\rm for}\
{\rm all}\ t\geq0\ {\rm a.s.}
\end{eqnarray*}

\end{thm}

We will give the proof of this theorem in the last section. In order
to find the stationary weak solution of SPDE (\ref{zhao21}), we need
to assume reasonable conditions on $f$ and $g$ so that we can check
the conditions in Theorem \ref{qi888}. Indeed under the weak
Lipschitz and monotone conditions posed in \cite{zh-zh}, all the
conditions of this theorem can be verified. In this paper, we will
consider the following conditions:
\begin{description}
\item[(A.1).] Functions $f$
and $g$ are
$\mathscr{B}_{\mathbb{R}^{d}}\otimes\mathscr{B}_{\mathbb{R}^{1}}\otimes\mathscr{B}_{\mathbb{R}^{d}}$
measurable and there exist constants $M,M_j,C,C_j,\alpha_j\geq0$
with $\sum_{j=1}^\infty M_j<\infty$, $\sum_{j=1}^\infty C_j<\infty$
and $\sum_{j=1}^\infty\alpha_j<{1\over2}$ s.t. for any $Y\in
L_{\rho}^2({\mathbb{R}^{d}};{\mathbb{R}^{1}})$, $X_1, X_2, Z_1,
Z_2\in L_{\rho}^2({\mathbb{R}^{d}};{\mathbb{R}^{d}})$, measurable
$U:{\mathbb{R}^{d}}\rightarrow [0,1]$,
\begin{eqnarray*}
&&\int_{\mathbb{R}^d}U(x)|f(X_1(x),Y(x),Z_1(x))-f(X_2(x),Y(x),Z_2(x))|^2\rho^{-1}(x)dx\\
&\leq&\int_{\mathbb{R}^d}U(x)\big(M|X_1(x)-X_2(x)|^2+C|Z_1(x)-Z_2(x)|^2\big)\rho^{-1}(x)dx,\\
&&\int_{\mathbb{R}^d}U(x)|g_j(X_1(x),Y_1(x),Z_1(x))-g_j(X_2(x),Y_2(x),Z_2(x))|^2\rho^{-1}(x)dx\\
&\leq&\int_{\mathbb{R}^d}U(x)\big(M_j|X_1(x)-X_2(x)|^2+C_j|Y_1(x)-Y_2(x)|^2\nonumber\\
&&\ \ \ \ \ \ \ \ \ \ \ \ \
+\alpha_j|Z_1(x)-Z_2(x)|^2\big)\rho^{-1}(x)dx.
\end{eqnarray*}
\item[(A.2).] For $p\in(2,q-1)$,
$\int_{\mathbb{R}^d}\|g(x,0,0)\|_{\mathcal{L}^p_{U_0}(\mathbb{R}^{1})}^p\rho^{-1}(x)dx<\infty$.
\item[(A.3).] There exists a constant
$M_0\geq0$ s.t. for any $t\geq0$, $x,z\in\mathbb{R}^{d}$,
$y\in\mathbb{R}^{1}$,
$$|f(x,y,z)|\leq M_0(1+|y|+|z|).$$
\item[(A.4).] There exists a
constant $\mu>0$ with
$2\mu-{pK}-pC-{{p(p-1)}\over2}\sum_{j=1}^{\infty}{C_j}>0$ s.t. for
any $Y_1, Y_2\in L_{\rho}^2({\mathbb{R}^{d}};{\mathbb{R}^{1}})$,
$X,Z\in L_{\rho}^2({\mathbb{R}^{d}};{\mathbb{R}^{d}})$, measurable
$U:{\mathbb{R}^{d}}\rightarrow [0,1]$,
\begin{eqnarray*}
&&\int_{\mathbb{R}^d}U(x)\big(Y_1(x)-Y_2(x)\big)\big(f(X(x),Y_1(x),Z(x))\nonumber\\
&&\ \ \ \ \ \ \ \ \ \ \ \ \ \ \ \ \ \ \ \ \ \ \ \ \ \ \ \ \ \ \ \ \ -f(X(x),Y_2(x),Z(x))\big)\rho^{-1}(x)dx\\
&\leq&-\mu\int_{\mathbb{R}^d}U(x){|Y_1(x)-Y_2(x)|}^2\rho^{-1}(x)dx.
\end{eqnarray*}
\item[(A.5).] For any $x\in\mathbb{R}^{d}$, $(y,z)\rightarrow f(x,y,z)$ is
continuous.
\item[(A.6).] The functions $b\in C_{l,b}^2(\mathbb{R}^{d};\mathbb{R}^{1})$, $\sigma\in
C_{l,b}^3(\mathbb{R}^{d}\times\mathbb{R}^{d};\mathbb{R}^{1})$, and
for $p$ given in {\rm(A.2)}, the global Lipschitz constant $L$ for
$b$ and $\sigma$ satisfies $K-pL-{p(p-1)\over2}L^2>0$.
\end{description}
Note here we don't assume $f$ is Lipschitz in the variable $y$. We
will prove
\begin{thm}\label{qi043} Under Conditions {\rm(A.1)}--{\rm(A.6)}, Eq.(\ref{qi13}) has a unique solution
$(Y_{s}^{t,x},Z_{s}^{t,x})$. Moreover
$E[\sup_{s\geq0}\int_{\mathbb{R}^{d}}{\rm
e}^{{-{pK}}s}{{|{Y}_s^{t,x}|}^p}\rho^{-1}(x)dx]<\infty$.
\end{thm}
\begin{thm}\label{qi044} Under Conditions {\rm(A.1)}--{\rm(A.6)}, let $u(t,\cdot)\triangleq Y_{t}^{t,\cdot}$, where
$(Y_{\cdot}^{t,\cdot},Z_{\cdot}^{t,\cdot})$ is the solution of
Eq.(\ref{qi13}). Then for arbitrary $T$ and $t\in[0,T]$,
$u(t,\cdot)$ is a weak solution for Eq.(\ref{zhang685}). Moreover,
$u(t,\cdot)$ is {\rm a.s.} continuous w.r.t. $t$ in
$L_{\rho}^2(\mathbb{R}^d;\mathbb{R}^1)$.
\end{thm}

It is obvious that the conditions of Theorem \ref{qi888} are
satisfied from the conclusions of Theorem \ref{qi044}, so we obtain
the stationary weak solution of SPDE (\ref{zhao21}):
\begin{cor}\label{zz44} Under Conditions {\rm(A.1)}--{\rm(A.6)}, for arbitrary $T$ and $t\in[0,T]$, let $v(t,\cdot)\triangleq Y_{T-t}^{T-t,\cdot}$, where $(Y_{\cdot}^{t,\cdot},Z_{\cdot}^{t,\cdot})$ is the solution of Eq.(\ref{qi13}) with $\hat{B}_s={B}_{T-s}-{B}_T$ for all $s\geq0$. Then $v(t,\cdot)$ is
a ``perfect" stationary weak solution of Eq.(\ref{zhao21}).
\end{cor}

In order to study the infinite horizon BDSDEs and stationary
solution of SPDEs, first we have to study the finite time BDSDEs and
therefore obtain the probabilistic representation of the weak
solutions of corresponding SPDEs. This will be given in Section
\ref{s32}. These results are novel, not only because of the
connection of BDSDEs and SPDEs, but also due to the fact that the
SPDEs we study appear to be new as coefficient $g$ of the noise can
be a general one. The existence and uniqueness of such equations
when $g$ is independent of $\nabla v$ or linearly dependent on
$\nabla v$ were studied by the pioneering works of Da Prato and
Zabczyk \cite{pr-za1}, Krylov \cite{krylov}. Our work shows that
studying the BDSDEs is a natural method for studying such general
SPDEs. The infinite horizon BDSDEs and stationary solution of SPDEs
will be studied in Section \ref{s38}.

\section{Finite Horizon BDSDEs and the Corresponding SPDEs}\label{s32}
\setcounter{equation}{0}

\subsection{Conditions and main results}\label{s33}
In this section, we will study the following BDSDEs on finite
horizon and establish its connection with SPDEs. This is necessary
to establish the solution of infinite horizon BDSDE and its
connection with the SPDEs.
\begin{eqnarray}\label{qi20}
Y_{s}^{t,x}&=&h(X_{T}^{t,x})+\int_{s}^{T}f(r,X_{r}^{t,x},Y_{r}^{t,x},Z_{r}^{t,x})dr\nonumber\\
&&-\int_{s}^{T}g(r,X_{r}^{t,x},Y_{r}^{t,x},Z_{r}^{t,x})d^\dagger{\hat{B}}_r-\int_{s}^{T}\langle
Z_{r}^{t,x},dW_r\rangle,\ \ \ 0\leq s\leq T.
\end{eqnarray}
Here $h:\Omega\times\mathbb{R}^{d}\longrightarrow{\mathbb{R}^1}$,
$f:[0,T]\times\mathbb{R}^{d}\times\mathbb{R}^1\times\mathbb{R}^{d}{\longrightarrow{\mathbb{R}^1}}$,
$g:[0,T]\times\mathbb{R}^{d}\times\mathbb{R}^1\times\mathbb{R}^{d}\longrightarrow
{\mathcal{L}^2_{U_0}(\mathbb{R}^{1})}$. Set $g_j\triangleq
g\sqrt{\lambda_j}e_j:[0,T]\times\mathbb{R}^{d}\times\mathbb{R}^1\times\mathbb{R}^{d}{\longrightarrow{\mathbb{R}^1}}$,
then Eq.(\ref{qi20}) is equivalent to
\begin{eqnarray*}\label{qi21}
Y_{s}^{t,x}&=&h(X_{T}^{t,x})+\int_{s}^{T}f(r,X_{r}^{t,x},Y_{r}^{t,x},Z_{r}^{t,x})dr\nonumber\\
&&-\sum_{j=1}^{\infty}\int_{s}^{T}g_j(r,X_{r}^{t,x},Y_{r}^{t,x},Z_{r}^{t,x})d^\dagger{\hat{\beta}}_j(r)-\int_{s}^{T}\langle
Z_{r}^{t,x},dW_r\rangle,\ \ \ 0\leq s\leq T.
\end{eqnarray*}
\begin{defi}\label{qi051}
A pair of processes $(Y_{\cdot}^{t,\cdot},Z_{\cdot}^{t,\cdot})\in
S^{2,0}([0,T];L_{\rho}^2({\mathbb{R}^{d}};{\mathbb{R}^{1}}))\bigotimes
M^{2,0}\\([0,T];L_{\rho}^2({\mathbb{R}^{d}};{\mathbb{R}^{d}}))$ is
called a solution of Eq.(\ref{qi20}) if for any $\varphi\in
C_c^{0}(\mathbb{R}^d;\mathbb{R}^1)$,
\begin{eqnarray}\label{qi22}
\int_{\mathbb{R}^{d}}Y_s^{t,x}\varphi(x)dx&=&\int_{\mathbb{R}^{d}}h(X_{T}^{t,x})\varphi(x)dx+\int_{s}^{T}\int_{\mathbb{R}^{d}}f(r,X_{r}^{t,x},Y_r^{t,x},Z_r^{t,x})\varphi(x)dxdr\nonumber\\
&&-\sum_{j=1}^{\infty}\int_{s}^{T}\int_{\mathbb{R}^{d}}g_j(r,X_{r}^{t,x},Y_r^{t,x},Z_r^{t,x})\varphi(x)dxd^\dagger{\hat{\beta}}_j(r)\nonumber\\
&&-\int_{s}^{T}\langle
\int_{\mathbb{R}^{d}}Z_r^{t,x}\varphi(x)dx,dW_r\rangle\ \ \ P-{\rm
a.s.}
\end{eqnarray}
\end{defi}
We assume
\begin{description}
\item[(H.1).] Function $h$ is $\mathscr{F}_{T,\infty}^{\hat{B}}\otimes\mathscr{B}_{\mathbb{R}^{d}}$ measurable and
$E[\int_{\mathbb{R}^{d}}|h(x)|^2\rho^{-1}(x)dx]<\infty$.
\item[(H.2).] Functions $f$ and $g$ are $\mathscr{B}_{[0,T]}\otimes\mathscr{B}_{\mathbb{R}^{d}}\otimes\mathscr{B}_{\mathbb{R}^{1}}\otimes\mathscr{B}_{\mathbb{R}^{d}}$ measurable and there exist constants
$C,C_j,\alpha_j\geq0$ with $\sum_{j=1}^\infty C_j<\infty$ and
$\sum_{j=1}^\infty\alpha_j<{1\over2}$ s.t. for any $r\in[0,T]$, $Y,
Y_1,Y_2\in L_{\rho}^2({\mathbb{R}^{d}};{\mathbb{R}^{1}})$,
$X,Z_1,Z_2\in L_{\rho}^2({\mathbb{R}^{d}};{\mathbb{R}^{d}})$,
\begin{eqnarray*}
&&\int_{\mathbb{R}^{d}}|f(r,X(x),Y(x),Z_1(x))-f(r,X(x),Y(x),Z_2(x))|^2\rho^{-1}(x)dx\\
&\leq&C\int_{\mathbb{R}^{d}}|Z_1(x)-Z_2(x)|^2\rho^{-1}(x)dx,\\
&&\int_{\mathbb{R}^{d}}|g_j(r,X(x),Y_1(x),Z_1(x))-g_j(r,X(x),Y_2(x),Z_2(x))|^2\rho^{-1}(x)dx\\
&\leq&\int_{\mathbb{R}^{d}}(C_j|Y_1(x)-Y_2(x)|^2+{\alpha_j}|Z_1(x)-Z_2(x)|^2)\rho^{-1}(x)dx.
\end{eqnarray*}
\item[(H.3).]
The integral
$\int_{0}^{T}\int_{\mathbb{R}^{d}}\|g(r,x,0,0)\|^2_{\mathcal{L}^2_{U_0}(\mathbb{R}^{1})}\rho^{-1}(x)dxdr<\infty$.
\item[(H.4).] There exists a constant
$M_0^{'}\geq0$ s.t. for any $r\in[0,T]$, $x,z\in\mathbb{R}^{d}$,
$y\in\mathbb{R}^{1}$,
$$|f(r,x,y,z)|\leq M_0^{'}(1+|y|+|z|).$$
\item[(H.5).] There exists a
constant $\mu\in\mathbb{R}^{1}$ s.t. for any $r\in[0,T]$, $Y_1,
Y_2\in L_{\rho}^2({\mathbb{R}^{d}};{\mathbb{R}^{1}})$, $X,Z\in
L_{\rho}^2({\mathbb{R}^{d}};{\mathbb{R}^{d}})$, measurable
$U:{\mathbb{R}^{d}}\rightarrow [0,1]$,
\begin{eqnarray*}
&&\int_{\mathbb{R}^d}U(x)\big(Y_1(x)-Y_2(x)\big)\big(f(r,X(x),Y_1(x),Z(x))\nonumber\\
&&\ \ \ \ \ \ \ \ \ \ \ \ \ \ \ \ \ \ \ \ \ \ \ \ \ \ \ \ \ \ \ \ \ -f(r,X(x),Y_2(x),Z(x))\big)\rho^{-1}(x)dx\\
&\leq&\mu\int_{\mathbb{R}^d}U(x){|Y_1(x)-Y_2(x)|}^2\rho^{-1}(x)dx.
\end{eqnarray*}
\item[(H.6).] For any $r\in[0,T]$, $x\in\mathbb{R}^{d}$, $(y,z)\rightarrow f(r,x,y,z)$ is
continuous.
\item[(H.7).] The functions $b\in C_{l,b}^2(\mathbb{R}^{d};\mathbb{R}^{d})$, $\sigma\in
C_{l,b}^3(\mathbb{R}^{d};\mathbb{R}^{d}\times\mathbb{R}^{d})$.
\end{description}

The first objective of this section is to prove
\begin{thm}\label{zz48} Under Conditions {\rm(H.1)}--{\rm(H.7)}, BDSDE (\ref{qi20}) has a unique solution.
\end{thm}
Then we will make the connection between the solutions of BDSDE
(\ref{qi20}) and SPDE (\ref{zhang685}).
\begin{thm}\label{zz45} Under Conditions {\rm(H.1)}--{\rm(H.7)},
if we define $u(t,x)=Y_t^{t,x}$, where $(Y_s^{t,x},Z_s^{t,x})$ is
the solution of Eq.(\ref{qi20}), then $u(t,x)$ is the unique weak
solution of Eq.(\ref{zhang685}) with $u(T,x)=h(x)$. Moreover,
$u(s,X_s^{t,x})=Y_s^{t,x}$, $(\sigma^*\nabla
u)(s,X^{t,x}_s)=Z_s^{t,x}$ for a.a. $s\in[t,T]$,
$x\in\mathbb{R}^{d}$ a.s.\\
\end{thm}

\subsection{Existence and uniqueness of the solutions of BDSDEs with finite dimensional noise}\label{s35}
In their pioneering work \cite{pa-pe1}, Pardoux and Peng solved the
following BSDE with Lipschitz conditions on the coefficient:
\begin{eqnarray}\label{zz54}
Y_{s}&=&\xi+\int_{s}^{T}f(r,Y_{r},Z_{r})dr-\int_{s}^{T}\langle
Z_{r},dW_r\rangle.
\end{eqnarray}
After that, many researchers studied how to weaken the Lipschitz
conditions so that the BSDE system can include more equations. To
name but a few, in \cite{pe}, \cite{le-sa}, \cite{pa}, \cite{ko},
\cite{br-de-hu-pr-st} and \cite{br-hu}, researchers made their
significant contributions to this subject. In \cite{le-sa},
Lepeltier and San Martin assumed that the $\mathbb{R}^{1}$-valued
function $f(r,y,z)$ satisfies the measurable condition, the $y$, $z$
linear growth condition and the $y$, $z$ continuous condition, then
they proved the existence of the solution of Eq.(\ref{zz54}). But
the uniqueness of solution failed to be proved since the comparison
theorem cannot be used under non-Lipschitz condition.

In \cite{shi}, after proving the comparison theorem of BDSDE with
Lipschitz condition, the authors used the method in \cite{le-sa} and
proved the corresponding result for the following
$\mathbb{R}^{1}$-valued BDSDE: 
\begin{eqnarray}\label{zz55}
Y_s=\xi+\int_{s}^{T}f(r,Y_r,Z_r)dr-\int_{s}^{T}\langle
g(r,Y_r,Z_r),d^\dagger\hat{B}_r\rangle-\int_{s}^{T}\langle
Z_r,dW_r\rangle.\ \
\end{eqnarray}
They assumed the same condition for $f$ as in \cite{le-sa} and
$g(r,y,z)$ satisfies the standard measurable condition and Lipschitz
condition w.r.t. $y$ and $z$. Then in Theorem 4.1 in \cite{shi},
they proved the existence of solution of Eq.(\ref{zz55}).

First we study the following BDSDE with finite dimensional noise
under non-Lipschitz conditions:
\begin{eqnarray}\label{zhang66100}
Y_s^{t,x,n}=&&h(X_{T}^{t,x})+\int_{s}^{T}f(r,X_{r}^{t,x},Y^{t,x,n}_r,Z^{t,x,n}_r)dr\\
&&-\sum_{j=1}^{n}\int_{s}^{T}g_j(r,X_{r}^{t,x},Y^{t,x,n}_r,Z^{t,x,n}_r)d^\dagger{\hat{\beta}}_j(r)-\int_{s}^{T}\langle
Z^{t,x,n}_r,dW_r\rangle.\nonumber
\end{eqnarray}
Note here in \cite{shi} and \cite{le-sa}, the authors only dealt
with the solution of Eq.(\ref{zhang66100}) for a fixed $x$ almost
surely. Of course if one is interested in the classical solution of
this SPDEs, it is easy to see that this implies one can solve
Eq.(\ref{zhang66100}) for all $x\in\mathbb{R}^{d}$ a.s. by some
standard perfection arguments. But we consider the solution in the
space
$S^{2,0}([0,T];L_{\rho}^2({\mathbb{R}^{d}};{\mathbb{R}^{1}}))\bigotimes
M^{2,0}([0,T];L_{\rho}^2({\mathbb{R}^{d}};{\mathbb{R}^{d}}))$ in
order to consider the weak solution of the SPDEs. The main task of
this subsection is to prove
\begin{thm}\label{zz000} Under Conditions {\rm(H.1)}--{\rm(H.7)},
Eq.(\ref{zhang66100}) has a unique solution
\begin{center}
$(Y_\cdot^{t,\cdot,n},Z_\cdot^{t,\cdot,n})\in
S^{2,0}([0,T];L_{\rho}^2({\mathbb{R}^{d}};{\mathbb{R}^{1}}))\bigotimes
M^{2,0}([0,T];L_{\rho}^2({\mathbb{R}^{d}};{\mathbb{R}^{d}}))$.
\end{center}
\end{thm}

We will first acknowledge the following Proposition \ref{zz001} at
the moment, then we prove Theorem \ref{zz000} with the help of
Proposition \ref{zz001}. Note that in the proof of Theorem
\ref{zz000} and Proposition \ref{zz001}, we can consider the
solution in
$S^{2,0}([t,T];L_{\rho}^2({\mathbb{R}^{d}};{\mathbb{R}^{1}}))\bigotimes
M^{2,0}([t,T];L_{\rho}^2({\mathbb{R}^{d}};{\mathbb{R}^{d}}))$ due to
the arguments in Remark 3.7 in \cite{zh-zh}.
\begin{prop}\label{zz001}
Given $(U_{\cdot}(\cdot),V_{\cdot}(\cdot))\in
S^{2,0}([0,T];L_{\rho}^2({\mathbb{R}^{d}};{\mathbb{R}^{1}}))\bigotimes
M^{2,0}([0,T];\\L_{\rho}^2({\mathbb{R}^{d}};{\mathbb{R}^{d}}))$,
then under Conditions {\rm(H.1)}--{\rm(H.7)}, the equation
\begin{eqnarray}\label{zz0}
Y_s^{t,x,n}=&&h(X_{T}^{t,x})+\int_{s}^{T}f(r,X_{r}^{t,x},Y^{t,x,n}_r,Z^{t,x,n}_r)dr\\
&&-\sum_{j=1}^{n}\int_{s}^{T}g_j(r,X_{r}^{t,x},U_r(x),V_r(x))d^\dagger{\hat{\beta}}_j(r)-\int_{s}^{T}\langle
Z^{t,x,n}_r,dW_r\rangle\nonumber
\end{eqnarray}
has a unique solution.
\end{prop}

{\em Proof of Theorem \ref{zz000}}. \underline{Uniqueness}. Assume
there exists another
$(\hat{Y}_{\cdot}^{t,\cdot,n},\hat{Z}_{\cdot}^{t,\cdot,n})\\\in
S^{2,0}([t,T];L_{\rho}^2({\mathbb{R}^{d}};{\mathbb{R}^{1}}))\bigotimes
M^{2,0}([t,T];L_{\rho}^2({\mathbb{R}^{d}};{\mathbb{R}^{d}}))$
satisfying (\ref{zhang66100}). Define
\begin{eqnarray*}
\bar{Y}_s^{t,x,n}=Y_{s}^{t,x,n}-\hat{Y}_{s}^{t,x,n}\ {\rm and}\
\bar{Z}_s^{t,x,n}=Z_{s}^{t,x,n}-\hat{Z}_{s}^{t,x,n},\ t\leq s\leq T.
\end{eqnarray*}
Then with probability $1$ we have that for a.e.
$x\in{\mathbb{R}^{d}}$, $(\bar{Y}_s^{t,x,n},\bar{Z}_s^{t,x,n})$
satisfies
\begin{eqnarray*}
\bar{Y}_s^{t,x,n}=&&\int_{s}^{T}\big(f(r,X_{r}^{t,x},Y^{t,x,n}_r,Z^{t,x,n}_r)-f(r,X_{r}^{t,x},\hat{Y}_{s}^{t,x,n},\hat{Z}_{s}^{t,x,n})\big)dr\nonumber\\
&&-\sum_{j=1}^{n}\int_{s}^{T}\big(g_j(r,X_{r}^{t,x},Y^{t,x,n}_r,Z^{t,x,n}_r)-g_j(r,X_{r}^{t,x},\hat{Y}_{s}^{t,x,n},\hat{Z}_{s}^{t,x,n})\big)d^\dagger{\hat{\beta}}_j(r)\nonumber\\
&&-\int_{s}^{T}\langle\bar{Z}^{t,x,n}_r,dW_r\rangle.\nonumber
\end{eqnarray*}
From Condition (H.4) and
$(\hat{Y}^{t,\cdot,n}_{\cdot},\hat{Z}^{t,\cdot,n}_{\cdot}),(Y^{t,\cdot,n}_{\cdot},Z^{t,\cdot,n}_{\cdot})\in
S^{2,0}([t,T];
L_{\rho}^2({\mathbb{R}^{d}};{\mathbb{R}^{1}}))\\\bigotimes
M^{2,0}([t,T];L_{\rho}^2({\mathbb{R}^{d}};{\mathbb{R}^{d}}))$, it
follows that
\begin{eqnarray*}
&&E[\int_{t}^{T}\int_{\mathbb{R}^{d}}|f(r,X_{r}^{t,x},Y^{t,x,n}_r,Z^{t,x,n}_r)-f(r,X_{r}^{t,x},\hat{Y}_{s}^{t,x,n},\hat{Z}_{s}^{t,x,n})|^2\rho^{-1}(x)dxdr]\nonumber\\
&\leq&C_pE[\int_{t}^{T}\int_{\mathbb{R}^{d}}(1+|Y^{t,x,n}_r|^2+|\hat{Y}^{t,x,n}_r|^2+|Z^{t,x,n}_r|^2+|\hat{Z}^{t,x,n}_r|^2)\rho^{-1}(x)dxdr]\nonumber\\
&<&\infty,
\end{eqnarray*}
where and in the rest of this paper $C_p$ is a generic constant. So
from Fubini theorem we have for a.e. $x\in{\mathbb{R}^{d}}$,
\begin{eqnarray*}
E[\int_{t}^{T}|f(r,X_{r}^{t,x},{Y}_{r}^{t,x,n},{Z}_{r}^{t,x,n})-f(r,X_{r}^{t,x},{\hat{Y}}_{r}^{t,x,n},{\hat{Z}}_{r}^{t,x,n})|^2dr]<\infty.
\end{eqnarray*}
Similarly, with Condition (H.2), we have for a.e.
$x\in{\mathbb{R}^{d}}$,
\begin{eqnarray*}
\sum_{j=1}^{n}E[\int_{t}^{T}|g_j(r,X_{r}^{t,x},{Y}_{r}^{t,x,n},{Z}_{r}^{t,x,n})-g_j(r,X_{r}^{t,x},{\hat{Y}}_{r}^{t,x,n},{\hat{Z}}_{r}^{t,x,n})|^2dr]<\infty.
\end{eqnarray*}
For a.e. $x\in{\mathbb{R}^{d}}$,
we apply the generalized It$\hat {\rm o}$'s formula
(\cite{el-tr-zh})
to ${\rm e}^{Ks}\psi_M\big(\bar{Y}_s^{t,x,n}\big)$, where
$K\in{\mathbb{R}^{1}}$ and
\begin{eqnarray*}
\psi_M(x)=x^2I_{\{-M\leq x<M\}}+M(2x-M)I_{\{x\geq
M\}}-M(2x+M)I_{\{x<-M\}}.
\end{eqnarray*}
Then
\begin{eqnarray}\label{zz1}
&&{\rm e}^{Ks}\psi_M(\bar{Y}_s^{t,x,n})+K\int_{s}^{T}{\rm e}^{Kr}\psi_M(\bar{Y}_r^{t,x,n})dr\nonumber\\
&&+\int_{s}^{T}{\rm e}^{Kr}I_{\{-M\leq\bar{Y}_r^{t,x,n}<M\}}|\bar{Z}_r^{t,x,n}|^2dr\nonumber\\
&=&\int_{s}^{T}{\rm e}^{Kr}\psi_M^{'}(\bar{Y}_r^{t,x,n})\big(f(r,X_{r}^{t,x},Y^{t,x,n}_r,Z^{t,x,n}_r)-f(r,X_{r}^{t,x},\hat{Y}_{r}^{t,x,n},\hat{Z}_{r}^{t,x,n})\big)dr\nonumber\\
&&+\sum_{j=1}^{n}\int_{s}^{T}{\rm e}^{Kr}I_{\{-M\leq\bar{Y}_r^{t,x,n}<M\}}|g_j(r,X_{r}^{t,x},Y^{t,x,n}_r,Z^{t,x,n}_r)\nonumber\\
&&\ \ \ \ \ \ \ \ \ \ \ \ \ \ \ \ \ \ \ \ \ \ \ \ \ \ \ \ \ \ \ \ \ \ \ \ \ \ \ \ \ \ -g_j(r,X_{r}^{t,x},\hat{Y}_{r}^{t,x,n},\hat{Z}_{r}^{t,x,n})|^2dr\nonumber\\
&&-\sum_{j=1}^{n}\int_{s}^{T}{\rm e}^{Kr}\psi_M^{'}(\bar{Y}_r^{t,x,n})\big(g_j(r,X_{r}^{t,x},Y^{t,x,n}_r,Z^{t,x,n}_r)\nonumber\\
&&\ \ \ \ \ \ \ \ \ \ \ \ \ \ \ \ \ \ \ \ \ \ \ \ \ \ \ \ \ \ \ \ \ \ -g_j(r,X_{r}^{t,x},\hat{Y}_{r}^{t,x,n},\hat{Z}_{r}^{t,x,n})\big)d^\dagger{\hat{\beta}}_j(r)\nonumber\\
&&-\int_{s}^{T}\langle{\rm
e}^{Kr}\psi_M^{'}(\bar{Y}_r^{t,x,n})\bar{Z}_r^{t,x,n},dW_r\rangle.
\end{eqnarray}
We can use the Fubini theorem to perfect (\ref{zz1}) so that
(\ref{zz1}) is satisfied for a.e. $x\in \mathbb{R}^d$ on a full
measure set that is independent of $x$. If we define
${{\psi_M^{'}(x)}\over x}=2$ when $x=0$, then
$0\leq{{\psi_M^{'}(\bar{Y}_r^{t,x,n})}\over{\bar{Y}_r^{t,x,n}}}\leq2$.
Taking integration over ${\mathbb{R}^{d}}$ on both sides of
(\ref{zz1}), we can apply the stochastic Fubini theorem
(\cite{pr-za1}).
Noting that the stochastic integrals are martingales, so taking the
expectation, we have
\begin{eqnarray*}\label{zhao666}
&&E[\int_{\mathbb{R}^{d}}{\rm e}^{Ks}\psi_M(\bar{Y}_s^{t,x,n})\rho^{-1}(x)dx]+KE[\int_{s}^{T}\int_{\mathbb{R}^{d}}{\rm e}^{Kr}\psi_M(\bar{Y}_r^{t,x,n})\rho^{-1}(x)dxdr]\nonumber\\
&&+E[\int_{s}^{T}\int_{\mathbb{R}^{d}}{\rm e}^{Kr}I_{\{-M\leq\bar{Y}_r^{t,x,n}<M\}}|\bar{Z}_r^{t,x,n}|^2\rho^{-1}(x)dxdr]\nonumber\\
&\leq&E[\int_{s}^{T}\int_{\mathbb{R}^{d}}{\rm e}^{Kr}{{\psi_M^{'}(\bar{Y}_r^{t,x,n})}\over{\bar{Y}_r^{t,x,n}}}\mu|\bar{Y}_r^{t,x,n}|^2\rho^{-1}(x)dxdr]\nonumber\\
&&+(2C+\sum_{j=1}^{\infty}C_j)E[\int_{s}^{T}\int_{\mathbb{R}^{d}}{\rm e}^{Kr}|\bar{Y}_{r}^{t,x,n}|^2\rho^{-1}(x)dxdr]\nonumber\\
&&+({1\over2}+\sum_{j=1}^{\infty}\alpha_j)E[\int_{s}^{T}\int_{\mathbb{R}^{d}}{\rm
e}^{Kr}|\bar{Z}_{r}^{t,x,n}|^2\rho^{-1}(x)dxdr].
\end{eqnarray*}
Taking the limit as $M\rightarrow\infty$ and applying the monotone
convergence theorem, we have
\begin{eqnarray}\label{zz2}
&&E[\int_{\mathbb{R}^{d}}{\rm e}^{Ks}|\bar{Y}_{s}^{t,x,n}|^2\rho^{-1}(x)dx]\nonumber\\
&&+({1\over2}-\sum_{j=1}^{\infty}\alpha_j)E[\int_{s}^{T}\int_{\mathbb{R}^{d}}{\rm e}^{Kr}|\bar{Z}_{r}^{t,x,n}|^2\rho^{-1}(x)dxdr]\nonumber\\
&&+(K-2\mu-2C-\sum_{j=1}^{\infty}C_j)E[\int_{s}^{T}\int_{\mathbb{R}^{d}}{\rm
e}^{Kr}|\bar{Y}_{r}^{t,x,n}|^2\rho^{-1}(x)dxdr]\leq0.
\end{eqnarray}
Note that all the terms on the left hand side of (\ref{zz2}) are
positive when $K$ is taken sufficiently large. So
by a ``standard" argument,
we have $\bar{Y}_s^{t,x,n}=0$ for all $s\in[t,T]$, a.a.
$x\in\mathbb{R}^{d}$ a.s. Also by (\ref{zz2}), for a.e.
$\bar{Z}_s^{t,x,n}=0$ for a.a. $x\in\mathbb{R}^{d}$ a.s. We can a.s.
We can modify the values of $Z$ at the measure zero exceptional set
of $s$ such that $\bar{Z}_s^{t,x,n}=0$ for all $x\in\mathbb{R}^{d}$
a.s.
\\

\underline{Existence}. If we regard Eq.(\ref{zz0}) as a mapping,
then by Proposition \ref{zz001}, (\ref{zz0}) is an iterated mapping
from
$S^{2,0}([t,T];L_{\rho}^2({\mathbb{R}^{d}};{\mathbb{R}^{1}}))\bigotimes
M^{2,0}([t,T];L_{\rho}^2({\mathbb{R}^{d}};{\mathbb{R}^{d}}))$ to
itself and we can obtain a sequence
$\{({Y}_r^{t,x,n,i},{Z}_r^{t,x,n,i})\}_{i=1}^\infty$ from this
mapping. We will prove that (\ref{zz0}) is a contraction mapping.
For this, define for $t\leq s\leq T$ and arbitrary given
$({Y}_\cdot^{t,\cdot,n,1},{Z}_\cdot^{t,\cdot,n,1})\in
S^{2,0}([t,T];L_{\rho}^2({\mathbb{R}^{d}};{\mathbb{R}^{1}}))\bigotimes
M^{2,0}\\([t,T];L_{\rho}^2({\mathbb{R}^{d}};{\mathbb{R}^{d}}))$,
\begin{eqnarray*}
&&\bar{Y}_s^{t,x,n,i}={Y}_s^{t,x,n,i}-{Y}_s^{t,x,n,i-1}, \ \
\bar{Z}_s^{t,x,n,i}={Z}_s^{t,x,n,i}-{Z}_s^{t,x,n,i-1},\\
&&\bar{g}^{i}_j(s,x)=g_j(s,X_{s}^{t,x},{Y}_s^{t,x,n,i},{Z}_s^{t,x,n,i})-g_j(s,X_{s}^{t,x},{Y}_s^{t,x,n,i-1},{Z}_s^{t,x,n,i-1}),\
i=2,3,\cdot\cdot\cdot
\end{eqnarray*}
Then for a.e. $x\in{\mathbb{R}^{d}}$,
$(\bar{Y}_s^{t,x,n,N},\bar{Z}_s^{t,x,n,N})$ satisfies
\begin{eqnarray*}
\bar{Y}_{s}^{t,x,n,N}&=&\int_{s}^{T}\big(f(r,X_{r}^{t,x},Y^{t,x,n,N}_r,Z^{t,x,n,N}_r)-f(r,X_{r}^{t,x},Y_{s}^{t,x,n,N-1},Z_{s}^{t,x,n,N-1})\big)dr\nonumber\\
&&-\sum_{j=1}^{n}\int_{s}^{T}\bar{g}^{N-1}_j(r,x)d^\dagger{\hat{\beta}}_j(r)-\int_{s}^{T}\langle\bar{Z}_{r}^{t,x,n,N},dW_r\rangle.
\end{eqnarray*}
Applying generalized the It$\hat {\rm o}$'s formula to ${\rm
e}^{Kr}\psi_M(\bar{Y}_{r}^{t,x,n,N})$ for a.e. $x\in\mathbb{R}^{d}$,
by the Young inequality, Condition (H.2) and (H.5), we can deduce
that
\begin{eqnarray}\label{zz4}
&&\int_{\mathbb{R}^{d}}{\rm e}^{Ks}\psi_M(\bar{Y}_s^{t,x,n,N})\rho^{-1}(x)dx+K\int_{s}^{T}\int_{\mathbb{R}^{d}}{\rm e}^{Kr}\psi_M(\bar{Y}_r^{t,x,n,N})\rho^{-1}(x)dxdr\nonumber\\
&&+\int_{s}^{T}\int_{\mathbb{R}^{d}}{\rm e}^{Kr}I_{\{-M\leq\bar{Y}_r^{t,x,n,N}<M\}}|\bar{Z}_r^{t,x,n,N}|^2\rho^{-1}(x)dxdr\nonumber\\
&\leq&\int_{s}^{T}\int_{\mathbb{R}^{d}}{\rm e}^{Kr}{{\psi_M^{'}(\bar{Y}_r^{t,x,n,N})}\over{\bar{Y}_r^{t,x,n,N}}}\mu|\bar{Y}_r^{t,x,n,N}|^2\rho^{-1}(x)dxdr\nonumber\\
&&+2C\int_{s}^{T}\int_{\mathbb{R}^{d}}{\rm e}^{Kr}|\bar{Y}_{r}^{t,x,n,N}|^2\rho^{-1}(x)dxdr\nonumber\\
&&+{1\over2}\int_{s}^{T}\int_{\mathbb{R}^{d}}{\rm e}^{Kr}|\bar{Z}_{r}^{t,x,n,N}|^2\rho^{-1}(x)dxdr\nonumber\\
&&+\sum_{j=1}^{\infty}C_j\int_{s}^{T}\int_{\mathbb{R}^{d}}{\rm e}^{Kr}|\bar{Y}_{r}^{t,x,n,N-1}|^2\rho^{-1}(x)dxdr\nonumber\\
&&+\sum_{j=1}^{\infty}\alpha_j\int_{s}^{T}\int_{\mathbb{R}^{d}}{\rm e}^{Kr}|\bar{Z}_{r}^{t,x,n,N-1}|^2\rho^{-1}(x)dxdr\nonumber\\
&&-\sum_{j=1}^{n}\int_{s}^{T}\int_{\mathbb{R}^{d}}{\rm e}^{Kr}\psi_M^{'}(\bar{Y}_r^{t,x,n,N})\bar{g}^{N-1}_j(r,x)\rho^{-1}(x)dxd^\dagger{\hat{\beta}}_j(r)\nonumber\\
&&-\int_{s}^{T}\langle\int_{\mathbb{R}^{d}}{\rm
e}^{Kr}\psi_M^{'}(\bar{Y}_r^{t,x,n,N})\bar{Z}_r^{t,x,n,N}\rho^{-1}(x)dx,dW_r\rangle.
\end{eqnarray}
Then taking expectation and the limit as $M\to\infty$, we have
\begin{eqnarray*}
&&(K-2\mu-2C)E[\int_{s}^{T}\int_{\mathbb{R}^{d}}{\rm e}^{Kr}|{\bar{Y}}_{r}^{t,x,n,N}|^2\rho^{-1}(x)dxdr]\nonumber\\
&&+{1\over2}E[\int_{s}^{T}\int_{\mathbb{R}^{d}}{\rm e}^{Kr}|{\bar{Z}}_{r}^{t,x,n,N}|^2\rho^{-1}(x)dxdr]\\
&\leq&E[\int_{s}^{T}\int_{\mathbb{R}^{d}}{\rm
e}^{Kr}(\sum_{j=1}^{\infty}C_j|{\bar{Y}}_{r}^{t,x,n,N-1}|^2+\sum_{j=1}^{\infty}\alpha_j|{\bar{Z}}_{r}^{t,x,n,N-1}|^2)\rho^{-1}(x)dxdr].
\end{eqnarray*}
First assuming that
$\sum_{j=1}^{\infty}C_j,\sum_{j=1}^{\infty}\alpha_j>0$, we have
\begin{eqnarray*}
&&(2K-4\mu-4C)E[\int_{s}^{T}\int_{\mathbb{R}^{d}}{\rm e}^{Kr}|{\bar{Y}}_{r}^{t,x,n,N}|^2\rho^{-1}(x)dxdr]\nonumber\\
&&+E[\int_{s}^{T}\int_{\mathbb{R}^{d}}{\rm e}^{Kr}|{\bar{Z}}_{r}^{t,x,n,N}|^2\rho^{-1}(x)dxdr]\\
&\leq&2\sum_{j=1}^{\infty}\alpha_jE[\int_{s}^{T}\int_{\mathbb{R}^{d}}{\rm
e}^{Kr}({{\sum_{j=1}^{\infty}C_j}\over{\sum_{j=1}^{\infty}\alpha_j}}|{\bar{Y}}_{r}^{t,x,n,N-1}|^2+|{\bar{Z}}_{r}^{t,x,n,N-1}|^2)\rho^{-1}(x)dxdr].
\end{eqnarray*}
Letting
$K=2\mu+2C+{{\sum_{j=1}^{\infty}C_j}\over{2\sum_{j=1}^{\infty}\alpha_j}}$,
we have
\begin{eqnarray}\label{zz5}
&&E[\int_{s}^{T}\int_{\mathbb{R}^{d}}{\rm e}^{Kr}({{\sum_{j=1}^{\infty}C_j}\over{\sum_{j=1}^{\infty}\alpha_j}}|{\bar{Y}}_{r}^{t,x,n,N}|^2+|{\bar{Z}}_{r}^{t,x,n,N}|^2)\rho^{-1}(x)dxdr]\\
&\leq&2\sum_{j=1}^{\infty}\alpha_jE[\int_{s}^{T}\int_{\mathbb{R}^{d}}{\rm
e}^{Kr}({{\sum_{j=1}^{\infty}C_j}\over{\sum_{j=1}^{\infty}\alpha_j}}|{\bar{Y}}_{r}^{t,x,n,N-1}|^2+|{\bar{Z}}_{r}^{t,x,n,N-1}|^2)\rho^{-1}(x)dxdr].\nonumber
\end{eqnarray}
Note that $E[\int_{t}^{T}\int_{\mathbb{R}^{d}}{\rm
e}^{Kr}({{\sum_{j=1}^{\infty}C_j}\over{\sum_{j=1}^{\infty}\alpha_j}}|\cdot|^2+|\cdot|^2)\rho^{-1}(x)dxdr]$
is equivalent to
$E[\int_{t}^{T}\int_{\mathbb{R}^{d}}\big(|\cdot|^2+|\cdot|^2\big)\rho^{-1}(x)dxdr]$.
From the contraction principle, the mapping (\ref{zz0}) has a pair
of fixed point
$(Y_{\cdot}^{t,\cdot,n,\infty},Z_{\cdot}^{t,\cdot,n,\infty})$ that
is the limit of the Cauchy sequence
${\{(Y_{\cdot}^{t,\cdot,n,N},Z_{\cdot}^{t,\cdot,n,N})\}}_{N=1}^{\infty}$
in
$M^{2,0}([t,T];L_{\rho}^2({\mathbb{R}^{d}};{\mathbb{R}^{1}}))\bigotimes
M^{2,0}([t,T];L_{\rho}^2\\({\mathbb{R}^{d}};{\mathbb{R}^{d}}))$. We
then prove that $Y_{\cdot}^{t,\cdot,n,\infty}$ is also the limit of
$Y_{\cdot}^{t,\cdot,n,N}$ in
$S^{2,0}([t,T];\\L_{\rho}^2({\mathbb{R}^{d}};{\mathbb{R}^{1}}))$ as
$N\rightarrow\infty$. For this, we only need to prove that
${\{Y_{\cdot}^{t,\cdot,n,N}\}}_{N=1}^{\infty}$ is a Cauchy sequence
in $S^{2,0}([t,T];L_{\rho}^2({\mathbb{R}^{d}};{\mathbb{R}^{1}}))$.
For this, from (\ref{zz4}), by the B-D-G inequality, the
Cauchy-Schwartz inequality and the Young inequality, we have
\begin{eqnarray*}\label{zz7}
&&E[\sup_{t\leq s\leq T}\int_{\mathbb{R}^{d}}{\rm e}^{Ks}\psi_M(\bar{Y}_s^{t,x,n,N})\rho^{-1}(x)dx]\nonumber\\
&\leq&C_pE[\int_{t}^{T}\int_{\mathbb{R}^{d}}(|\bar{Y}_{r}^{t,x,n,N}|^2+|\bar{Z}_{r}^{t,x,n,N}|^2\nonumber\\
&&\ \ \ \ \ \ \ \ \ \ \ \ \ \ \ \ \ \ +|\bar{Y}_{r}^{t,x,n,N-1}|^2+|\bar{Z}_{r}^{t,x,n,N-1}|^2)\rho^{-1}(x)dxdr]\nonumber\\
&&+C_pE[\sqrt{\int_{t}^{T}\int_{\mathbb{R}^{d}}|\psi_M^{'}(\bar{Y}_r^{t,x,n,N})|^2\rho^{-1}(x)dx\int_{\mathbb{R}^{d}}\sum_{j=1}^{n}|\bar{g}^{N-1}_j(r,x)|^2\rho^{-1}(x)dxdr}]\nonumber\\
&&+C_pE[\sqrt{\int_{t}^{T}\int_{\mathbb{R}^{d}}|\psi_M^{'}(\bar{Y}_r^{t,x,n,N})|^2\rho^{-1}(x)dx\int_{\mathbb{R}^{d}}|\bar{Z}_r^{t,x,n,N}|^2\rho^{-1}(x)dxdr}]\nonumber\\
&\leq&C_pE[\int_{t}^{T}\int_{\mathbb{R}^{d}}(|\bar{Y}_{r}^{t,x,n,N}|^2+|\bar{Z}_{r}^{t,x,n,N}|^2\nonumber\\
&&\ \ \ \ \ \ \ \ \ \ \ \ \ \ \ \ \ \ +|\bar{Y}_{r}^{t,x,n,N-1}|^2+|\bar{Z}_{r}^{t,x,n,N-1}|^2)\rho^{-1}(x)dxdr]\nonumber\\
&&+{1\over5}E[\sup_{t\leq s\leq
T}\int_{\mathbb{R}^{d}}|\psi_M^{'}(Y_s(x))|^2\rho^{-1}(x)dx],
\end{eqnarray*}
where $C_p$ depends on $|\mu|$, $C$, $\sum_{j=1}^{\infty}\alpha_j$,
$\sum_{j=1}^{\infty}C_j$ and the fixed constant in the B-D-G
inequality. Taking the limit as $M\rightarrow\infty$ and applying
the monotone convergence theorem, we have
\begin{eqnarray}\label{zz6}
&&E[\sup_{t\leq s\leq T}\int_{\mathbb{R}^{d}}{\rm e}^{Ks}|\bar{Y}_{s}^{t,x,n,N}|^2\rho^{-1}(x)dx]\\
&\leq&M_0^{''}E[\int_{s}^{T}\int_{\mathbb{R}^{d}}{\rm
e}^{Kr}\big(|\bar{Y}_{r}^{t,x,n,N-1}|^2+|\bar{Z}_{r}^{t,x,n,N-1}|^2\nonumber\\
&&\ \ \ \ \ \ \ \ \ \ \ \ \ \ \ \ \ \ \ \ \ \ \ \ \
+|\bar{Y}_{r}^{t,x,n,N}|^2+|\bar{Z}_{r}^{t,x,n,N}|^2\big)\rho^{-1}(x)dxdr],\nonumber
\end{eqnarray}
where $M_0^{''}>0$ is independent of $n$ and $N$. Without losing any
generality, assume that $M\geq N$. We can deduce from (\ref{zz5})
and (\ref{zz6}) that
\begin{eqnarray*}
&&\big(E[\sup_{t\leq s\leq T}\int_{\mathbb{R}^{d}}{{|{Y}_{s}^{t,x,n,M}-{Y}_{s}^{t,x,n,N}|}^2}\rho^{-1}(x)dx]\big)^{1\over2}\nonumber\\
&\leq&\sum_{i=N+1}^{M}\big(E[\sup_{t\leq s\leq T}\int_{\mathbb{R}^{d}}{{|\bar{Y}_{s}^{t,x,n,i}|}^2}\rho^{-1}(x)dx]\big)^{1\over2}\nonumber\\
&\leq&\sum_{i=N+1}^{M}\big(M_0^{''}E[\int_{t}^{T}\int_{\mathbb{R}^{d}}{\rm e}^{Kr}\big(|\bar{Y}_{r}^{t,x,n,i-1}|^2+|\bar{Z}_{r}^{t,x,n,i-1}|^2\nonumber\\
&&\ \ \ \ \ \ \ \ \ \ \ \ \ \ \ \ \ \ \ \ \ \ \ \ \ \ \ \ \ \ \ \ \ \ \ +|\bar{Y}_{r}^{t,x,n,i}|^2+|\bar{Z}_{r}^{t,x,n,i}|^2\big)\rho^{-1}(x)dxdr]\big)^{1\over2}\nonumber\\
&\leq&\sum_{i=N+1}^{M}\big((1+{{\sum_{j=1}^{\infty}\alpha_j}\over{\sum_{j=1}^{\infty}C_j}})M_0^{''}E[\int_{t}^{T}\int_{\mathbb{R}^{d}}{\rm e}^{Kr}\big({{\sum_{j=1}^{\infty}C_j}\over{\sum_{j=1}^{\infty}\alpha_j}}|\bar{Y}_{r}^{t,x,n,i-1}|^2+|\bar{Z}_{r}^{t,x,n,i-1}|^2\nonumber\\
&&\ \ \ \ \ \ \ \ \ \ \ \ \ \ \ \ \ \ \ \ \ \ \ \ \ \ \ \ \ \ \ \ \ \ \ \ \ +{{\sum_{j=1}^{\infty}C_j}\over{\sum_{j=1}^{\infty}\alpha_j}}|\bar{Y}_{r}^{t,x,n,i}|^2+|\bar{Z}_{r}^{t,x,n,i}|^2\big)\rho^{-1}(x)dxdr]\big)^{1\over2}\nonumber\\
&\leq&\sum_{i=N+1}^{M}\big((2+{{2\sum_{j=1}^{\infty}\alpha_j}\over{\sum_{j=1}^{\infty}C_j}})M_0^{''}\nonumber\\
&&\times E[\int_{t}^{T}\int_{\mathbb{R}^{d}}{\rm e}^{Kr}({{\sum_{j=1}^{\infty}C_j}\over{\sum_{j=1}^{\infty}\alpha_j}}|\bar{Y}_{r}^{t,x,n,i-1}|^2+|\bar{Z}_{r}^{t,x,n,i-1}|^2)\rho^{-1}(x)dxdr]\big)^{1\over2}\nonumber\\
&\leq&\sum_{i=N+1}^{\infty}({2\sum_{j=1}^{\infty}\alpha_j})^{{i-2}\over2}\big((2+{{2\sum_{j=1}^{\infty}\alpha_j}\over{\sum_{j=1}^{\infty}C_j}})M_0^{''}\nonumber\\
&&\times E[\int_{t}^{T}\int_{\mathbb{R}^{d}}{\rm
e}^{Kr}({{\sum_{j=1}^{\infty}C_j}\over{\sum_{j=1}^{\infty}\alpha_j}}|{Y}_{r}^{t,x,n,1}|^2+|{Z}_{r}^{t,x,n,1}|^2)\rho^{-1}(x)dxdr]\big)^{1\over2}\longrightarrow0\nonumber
\end{eqnarray*}
as $M$, $N\longrightarrow\infty$, since
$2\sum_{j=1}^{\infty}\alpha_j<1$. So we proved our claim.

If either or both
$\sum_{j=1}^{\infty}C_j,\sum_{j=1}^{\infty}\alpha_j=0$, we can prove
the above convergence using similar method or the above convergence
is trivially correct. Theorem \ref{zz000} is proved. $\hfill\diamond$\\

The remaining work in this subsection is to prove Proposition
\ref{zz001}. First we do some preparations.
\begin{lem}\label{zz002} Under Conditions {\rm(H.1)}--{\rm(H.7)}, if there exists $(Y_\cdot(\cdot),Z_\cdot(\cdot))\in
M^{2,0}([t,T];L_{\rho}^2({\mathbb{R}^{d}};{\mathbb{R}^{1}}))\bigotimes
M^{2,0}([t,T];L_{\rho}^2({\mathbb{R}^{d}};{\mathbb{R}^{d}}))$
satisfying the spatial integral form of Eq.(\ref{zhang66100}) for
$t\leq s\leq T$, then $Y_\cdot(\cdot)\in
S^{2,0}([t,T];L_{\rho}^2({\mathbb{R}^{d}};{\mathbb{R}^{1}}))$ and
therefore $(Y_s(x),Z_s(x))$ is a solution of Eq.(\ref{zhang66100}).
\end{lem}
{\em Proof}. Similar to the proof of Lemma 3.3 in \cite{zh-zh}, we
can prove ${Y}_{s}(\cdot)$ is continuous w.r.t. $s$ in
$L_{\rho}^2({\mathbb{R}^{d}};{\mathbb{R}^{1}})$ under the conditions
of this lemma. We only mention that we can use Condition (H.4) to
deal with the term $f(r,X_{r}^{t,x},Y_r(x),Z_r(x))$ although there
is no weak Lipschitz condition for $Y_r(x)$. We omit the proof here.
Now we only show the proof of
$E[\sup_{t\leq s\leq
T}\int_{\mathbb{R}^{d}}|Y_s(x)|^2\rho^{-1}(x)dx]<\infty$ briefly.
For a.e. $x\in{\mathbb{R}^{d}}$, applying the generalized It$\hat
{\rm o}$'s formula to $\psi_M\big(Y_r(x)\big)$, by Lemma
\ref{qi045},
the B-D-G inequality and the Cauchy-Schwartz inequality, we have
\begin{eqnarray*}
&&E[\sup_{t\leq s\leq T}\int_{\mathbb{R}^{d}}\psi_M(Y_s(x))\rho^{-1}(x)dx]\nonumber\\
&\leq&C_pE[\int_{\mathbb{R}^{d}}{|h(x)|^2}\rho^{-1}(x)dx]+C_pE[\int_{t}^{T}\int_{\mathbb{R}^{d}}(|Y_r(x)|^2+|Z_r(x)|^2)\rho^{-1}(x)dxdr]\nonumber\\
&&+C_p\sum_{j=1}^{n}\int_{t}^{T}\int_{\mathbb{R}^{d}}(1+|g_j(r,x,0,0)|^2)\rho^{-1}(x)dxdr<\infty.
\end{eqnarray*}
So taking the limit as $M\rightarrow\infty$ and applying the
monotone convergence theorem, we have $Y_\cdot(\cdot)\in
S^{2,0}([t,T];L_{\rho}^2({\mathbb{R}^{d}};{\mathbb{R}^{1}}))$.
Recall that a solution of Eq.(\ref{zhang66100}) is a pair of
processes in
$S^{2,0}([0,T];L_{\rho}^2({\mathbb{R}^{d}};{\mathbb{R}^{1}}))\bigotimes
M^{2,0}([0,T];L_{\rho}^2({\mathbb{R}^{d}};{\mathbb{R}^{d}}))$
satisfying the spatial integral form of Eq.(\ref{zhang66100}),
therefore $(Y_s(x),Z_s(x))$ is a solution of
Eq.(\ref{zhang66100}). $\hfill\diamond$\\

From the proof of Lemma \ref{zz002}, one can similarly deduce that
\begin{cor}\label{zz003} Under Conditions {\rm(H.1)}--{\rm(H.7)}, if there exists
$(Y_\cdot(\cdot),Z_\cdot(\cdot))\in
M^{2,0}([t,T];L_{\rho}^2({\mathbb{R}^{d}};{\mathbb{R}^{1}}))\bigotimes
M^{2,0}([t,T];L_{\rho}^2({\mathbb{R}^{d}};{\mathbb{R}^{d}}))$
satisfying the spatial integral form of Eq.(\ref{zz0}) for $t\leq
s\leq T$, then $Y_\cdot(\cdot)\in
S^{2,0}([t,T];L_{\rho}^2({\mathbb{R}^{d}};{\mathbb{R}^{1}}))$ and
therefore $(Y_s(x),Z_s(x))$ is a solution of Eq.(\ref{zz0}).
\end{cor}

For the rest of this paper, we will leave out the similar
localization argument as in the proof of Theorem \ref{zz000} and
Lemma \ref{zz002} when applying It$\hat {\rm o}$'s formula to save
the space of this paper.
\\

{\em Proof of Proposition \ref{zz001}}. The proof of the uniqueness
is rather similar to the uniqueness proof in Theorem \ref{zz000}, so
it is omitted.

\underline{Existence}. Define
\begin{eqnarray*}
\tilde{f}^{x}(r,y,z)=f(r,X_{r}^{t,x},y,z)\ {\rm and}\
\tilde{g}_j^{x}(r)=g_j(r,X_{r}^{t,x},U_r(x),V_r(x)),
\end{eqnarray*}
then for a.e. $x\in\mathbb{R}^{d}$, (\ref{zz0}) becomes
\begin{eqnarray}\label{zz9}
Y_s^{t,x,n}=&&h(X_{T}^{t,x})+\int_{s}^{T}\tilde{f}^{x}(r,Y^{t,x,n}_r,Z^{t,x,n}_r)dr\nonumber\\
&&-\sum_{j=1}^{n}\int_{s}^{T}\tilde{g}_j^{x}(r)d^\dagger{\hat{\beta}}_j(r)-\int_{s}^{T}\langle
Z^{t,x,n}_r,dW_r\rangle.
\end{eqnarray}
Then it is easy to see that for a.e. $x\in\mathbb{R}^{d}$,
$\tilde{f}^{x}$ and $\tilde{g}_j^{x}$ satisfy
\begin{description}
\item[(H.1)$'$.] $\tilde{f}^{x}:[t,T]\times\Omega\times\mathbb{R}^1\times\mathbb{R}^d{\longrightarrow{\mathbb{R}^1}}$ is $\mathscr{B}_{[t,T]}\otimes\mathscr{F}_{s,T}\bigvee{\mathscr{F}_{T,\infty}^{\hat{B}}}\otimes\mathscr{B}_{\mathbb{R}^{1}}\otimes\mathscr{B}_{\mathbb{R}^{d}}$ measurable and\\
$\tilde{g}_j^{x}:[t,T]\times\Omega{\longrightarrow{\mathbb{R}^1}}$
is
$\mathscr{B}_{[t,T]}\otimes\mathscr{F}_{s,T}\bigvee{\mathscr{F}_{T,\infty}^{\hat{B}}}$
measurable.
\item[(H.2)$'$.] For any $r\in[t,T]$, $y\in\mathbb{R}^{1}$,
$|\tilde{f}^{x}(r,y,z)|\leq M_0^{'}(1+|y|+|z|)$.
\item[(H.3)$'$.] For any $r\in[t,T]$, $(y,z)\rightarrow\tilde{f}^{x}(r,y,z)$ is
continuous.
\end{description}
By Theorem 4.1 in \cite{shi}, for a.e. $x\in\mathbb{R}^{d}$,
Eq.(\ref{zz9}), as well as Eq.(\ref{zz0}), has a solution
$(Y_s^{t,x,n},Z_s^{t,x,n})\in
M^{2,0}([t,T];{\mathbb{R}^{1}})\bigotimes
M^{2,0}([t,T];{\mathbb{R}^{d}})$. In the following, we will prove
that $(Y_s^{t,x,n},Z_s^{t,x,n})\in
M^{2,0}([t,T];L_{\rho}^2({\mathbb{R}^{d}};{\mathbb{R}^{d}}))\bigotimes
M^{2,0}([t,T];L_{\rho}^2\\({\mathbb{R}^{d}};{\mathbb{R}^{1}}))$
under the conditions of Proposition \ref{zz001}.

First by Condition (H.4) or Condition (H.2)$'$, Conditions (H.2),
(H.3) and (H.7), for a.e. $x\in{\mathbb{R}^{d}}$, we have
\begin{eqnarray*}
&&E[\int_{t}^{T}|f(r,X_{r}^{t,x},Y^{t,x,n}_r,Z^{t,x,n}_r)|^2dr]\nonumber\\
&&+\sum_{j=1}^{n}E[\int_{t}^{T}|g_j(r,X_{r}^{t,x},U_r(x),V_r(x))|^2dr]<\infty.
\end{eqnarray*}
Then for a.e. $x\in{\mathbb{R}^{d}}$, applying the generalized
It$\hat {\rm o}$'s formula to ${\rm e}^{Kr}|Y_r^{t,x,n}|^2$, we have
\begin{eqnarray*}
&&E[{\rm e}^{Ks}|{Y}_s^{t,x,n}|^2]+KE[\int_{s}^{T}{\rm e}^{Kr}|{Y}_r^{t,x,n}|^2dr]+E[\int_{s}^{T}{\rm e}^{Kr}|{Z}_r^{t,x,n}|^2dr]\nonumber\\
&=&E[{\rm e}^{KT}|h(X_{T}^{t,x})|^2]+2E[\int_{s}^{T}{\rm e}^{Kr}{Y}_r^{t,x,n}f(r,X_{r}^{t,x},Y^{t,x,n}_r,Z^{t,x,n}_r)dr]\nonumber\\
&&+\sum_{j=1}^{n}E[\int_{s}^{T}{\rm
e}^{Kr}|g_j(r,X_{r}^{t,x},U_r(x),V_r(x))|^2dr].\nonumber
\end{eqnarray*}
Taking the integration over ${\mathbb{R}^{d}}$ and by Conditions
(H.1)--(H.5), (H.7) and Lemma \ref{qi045}, we have
\begin{eqnarray*}
&&E[\int_{\mathbb{R}^{d}}{\rm e}^{Ks}|{Y}_s^{t,x,n}|^2\rho^{-1}(x)dx]+KE[\int_{s}^{T}\int_{\mathbb{R}^{d}}{\rm e}^{Kr}|{Y}_r^{t,x,n}|^2\rho^{-1}(x)dxdr]\nonumber\\
&&+E[\int_{s}^{T}\int_{\mathbb{R}^{d}}{\rm e}^{Kr}|{Z}_r^{t,x,n}|^2\rho^{-1}(x)dxdr]\nonumber\\
&=&E[\int_{\mathbb{R}^{d}}{\rm e}^{KT}|h(X_{T}^{t,x})|^2\rho^{-1}(x)dx]\nonumber\\
&&+2E[\int_{s}^{T}\int_{\mathbb{R}^{d}}{\rm e}^{Kr}{Y}_r^{t,x,n}f(r,X_{r}^{t,x},Y^{t,x,n}_r,Z^{t,x,n}_r)\rho^{-1}(x)dxdr]\nonumber\\
&&+\sum_{j=1}^{n}E[\int_{s}^{T}\int_{\mathbb{R}^{d}}{\rm e}^{Kr}|g_j(r,X_{r}^{t,x},U_r(x),V_r(x))|^2\rho^{-1}(x)dxdr]\nonumber\\
&\leq&C_pE[\int_{\mathbb{R}^{d}}|h(x)|^2\rho^{-1}(x)dx]\nonumber\\
&&+(2\mu+2C+1)E[\int_{s}^{T}\int_{\mathbb{R}^{d}}{\rm e}^{Kr}|{Y}_r^{t,x,n}|^2\rho^{-1}(x)dxdr]\nonumber\\
&&+{1\over2}E[\int_{s}^{T}\int_{\mathbb{R}^{d}}{\rm e}^{Kr}|{Z}_r^{t,x,n}|^2\rho^{-1}(x)dxdr]+C_p\nonumber\\
&&+C_pE[\int_{s}^{T}\int_{\mathbb{R}^{d}}(|U_r(x)|^2+|V_{r}(x)|^2)\rho^{-1}(x)dxdr]\nonumber\\
&&+C_p\sum_{j=1}^{n}\int_{s}^{T}\int_{\mathbb{R}^{d}}|g_j(r,x,0,0)|^2\rho^{-1}(x)dxdr.\nonumber
\end{eqnarray*}
It turns out that
\begin{eqnarray}\label{zz11}
&&E[\int_{\mathbb{R}^{d}}{\rm e}^{Ks}|{Y}_s^{t,x,n}|^2\rho^{-1}(x)dx]\nonumber\\
&&+(K-2\mu-2C-1)E[\int_{s}^{T}\int_{\mathbb{R}^{d}}{\rm e}^{Kr}|{Y}_r^{t,x,n}|^2\rho^{-1}(x)dxdr]\nonumber\\
&&+{1\over2}E[\int_{s}^{T}\int_{\mathbb{R}^{d}}{\rm e}^{Kr}|{Z}_r^{t,x,n}|^2\rho^{-1}(x)dxdr]\nonumber\\
&\leq&C_pE[\int_{\mathbb{R}^{d}}|h(x)|^2\rho^{-1}(x)dx]+C_p\nonumber\\
&&+C_pE[\int_{s}^{T}\int_{\mathbb{R}^{d}}(|U_r(x)|^2+|V_{r}(x)|^2)\rho^{-1}(x)dxdr]\nonumber\\
&&+C_p\sum_{j=1}^{n}\int_{s}^{T}\int_{\mathbb{R}^{d}}|g_j(r,x,0,0)|^2\rho^{-1}(x)dxdr\nonumber\\
&<&\infty.
\end{eqnarray}
Taking $K$ sufficiently large, we can see that
$(Y_\cdot^{t,\cdot,n},Z_\cdot^{t,\cdot,n})\in
M^{2,0}([t,T];L_{\rho}^2({\mathbb{R}^{d}};\\{\mathbb{R}^{d}}))\bigotimes
M^{2,0}([t,T];L_{\rho}^2({\mathbb{R}^{d}};{\mathbb{R}^{1}}))$ and
for a.e. $x\in{\mathbb{R}^{d}}$, $(Y_s^{t,x,n},Z_s^{t,x,n})$
satisfies Eq.(\ref{zz0}) on a full measure set
$\Omega^x\subset\Omega$ dependent on $x$. But we can use the Fubini
theorem to perfect Eq.(\ref{zz0}) so that
$(Y_s^{t,x,n},Z_s^{t,x,n})$ satisfies (\ref{zz0}) for a.e. $x\in
\mathbb{R}^d$ on a full measure set $\tilde{\Omega}$ independent of
$x$. To see this, from (\ref{zz11}), we have for any $s\in[t,T]$,
\begin{eqnarray}\label{999}
E[\int_{\mathbb{R}^{d}}{\rm
e}^{Ks}|{Y}_s^{t,x,n}|^2\rho^{-1}(x)dx]=\int_{\mathbb{R}^{d}}E[{\rm
e}^{Ks}|{Y}_s^{t,x,n}|^2\rho^{-1}(x)]dx<\infty,\
\end{eqnarray}
so for a.e. $x\in{\mathbb{R}^{d}}$, there exists a full measure set
$\Omega^x\subset\Omega$ s.t. ${Y}_s^{t,x,n}<\infty$ on $\Omega^x$.
Denote the right hand side of (\ref{zz0}) by $F(s,x)$.
Then by Eq.(\ref{zz0}), for $x\in\mathbb{R}^{d}$, there exists a
full measure set ${\Omega'}^x\subset\Omega$ s.t.
$Y_s^{t,x,n}=F(s,x)$ on ${\Omega'}^x$. Then for a.e.
$x\in{\mathbb{R}^{d}}$, we have $Y_s^{t,x,n}=F(s,x)$ on
$\Omega^x\bigcap{\Omega'}^x$. Since now for a.e.
$x\in{\mathbb{R}^{d}}$, ${Y}_s^{t,x,n}<\infty$ on
$\Omega^x\bigcap{\Omega'}^x$, so $F(s,x)<\infty$ and we can move
$F(s,x)$ to the other side of the equality to have
$Y_s^{t,x,n}-F(s,x)=0$ on the full measure set
$\Omega^x\bigcap{\Omega'}^x$. Thus
\begin{eqnarray*}
\int_{\mathbb{R}^{d}}E[|Y_s^{t,x,n}-F(s,x)|]dx=0.
\end{eqnarray*}
By the Fubini theorem, we have
\begin{eqnarray*}
E[\int_{\mathbb{R}^{d}}|Y_s^{t,x,n}-F(s,x)|dx]=0.
\end{eqnarray*}
This means that there exists a full measure set $\dot{\Omega}$
independent of $x$ s.t. on $\dot{\Omega}$, $Y_s^{t,x,n}-F(s,x)=0$
for $x\in\dot{{\mathcal{E}}}^\omega$, where
$\dot{{\mathcal{E}}}^\omega$ is a full measure set in
$\mathbb{R}^{d}$ and depends on $\omega$. Similarly, from
(\ref{999}), we also know that there exists another full measure set
$\ddot{\Omega}$ independent of $x$ s.t. on $\ddot{\Omega}$,
$Y_s^{t,x,n}<\infty$ for $x\in\ddot{{\mathcal{E}}}^\omega$, where
$\ddot{{\mathcal{E}}}^\omega$ is a full measure set in
$\mathbb{R}^{d}$ and depends on $\omega$. Take
$\tilde{\Omega}=\dot{\Omega}\bigcap\ddot{\Omega}$ and
$\tilde{{\mathcal{E}}}^\omega=\dot{{\mathcal{E}}}^\omega\bigcap\ddot{{\mathcal{E}}}^\omega$,
then both are still a full measure set and on $\tilde{\Omega}$,
$Y_s^{t,x,n}<\infty$ for $x\in\tilde{{\mathcal{E}}}^\omega$,
furthermore $F(s,x)<\infty$. We can move items in the equality
$Y_s^{t,x,n}-F(s,x)=0$ to have $Y_s^{t,x,n}=F(s,x)$ for $x\in
\tilde{{\mathcal{E}}}^\omega$ on a full measure set $\tilde{\Omega}$
independent of $x$.

Now we have $(Y_s^{t,x,n},Z_s^{t,x,n})\in
M^{2,0}([t,T];L_{\rho}^2({\mathbb{R}^{d}};{\mathbb{R}^{d}}))\bigotimes
M^{2,0}([t,T];L_{\rho}^2\\({\mathbb{R}^{d}};{\mathbb{R}^{1}}))$ and
for $t\leq s\leq T$, $(Y_s^{t,x,n},Z_s^{t,x,n})$ satisfies
(\ref{zz0}) for a.e. $x\in \mathbb{R}^d$ on a full measure set
$\tilde{\Omega}$ independent of $x$. Then for any $\varphi\in
C_c^{0}(\mathbb{R}^d;\mathbb{R}^1)$, multiplying by $\varphi$ on
both sides of Eq.(\ref{zz0}) and taking the integration over
$\mathbb{R}^d$, we have $(Y_s^{t,x,n},Z_s^{t,x,n})$ satisfies the
spatial integral form of Eq.(\ref{zz0}) for $t\leq s\leq T$. By
Corollary \ref{zz003}, $Y_\cdot^{t,\cdot,n}\in
S^{2,0}([t,T];L_{\rho}^2({\mathbb{R}^{d}};{\mathbb{R}^{1}}))$ and
$(Y_s^{t,x,n},Z_s^{t,x,n})$ is a solution of Eq.(\ref{zz0}). $\hfill\diamond$\\

\subsection{Existence and uniqueness of solutions of BDSDEs with infinite dimensional noise}\label{s36}
Following a similar procedure as in the proof of Lemma \ref{zz002},
and applying It$\hat {\rm o}$'s formula to ${\rm
e}^{Kr}{|{Y}_{r}^{t,x,n}|}^2$, by the B-D-G inequality we have the
following estimation for the solution of Eq.(\ref{zhang66100}):
\begin{prop}\label{zz58}
Under the conditions of Theorem \ref{zz48},
$({Y}_s^{t,x,n},Z^{t,x,n}_{s})$ satisfies
\begin{eqnarray*}
\sup_nE[\sup_{0\leq s\leq
T}\int_{\mathbb{R}^d}|Y_s^{t,x,n}|^2\rho^{-1}(x)dx]+\sup_nE[\int_{0}^{T}\int_{\mathbb{R}^{d}}|Z_r^{t,x,n}|^2\rho^{-1}(x)dxdr]<\infty.
\end{eqnarray*}
\end{prop}
Now we turn to the proof of the first main theorem of this section.
\\

{\em Proof of Theorem \ref{zz48}}. The proof of the uniqueness is
rather similar to the uniqueness proof in Theorem \ref{zz000}, so it
is omitted.

\underline{Existence}. By Theorem \ref{zz000}, for each $n$, there
exists a unique solution
$({Y}_\cdot^{t,\cdot,n},{Z}_\cdot^{t,\cdot,n})$ to
Eq.(\ref{zhang66100}), so
$(Y^{t,\cdot,n}_{\cdot},Z^{t,\cdot,n}_{\cdot})\in S^{2,0}([0,T];
L_{\rho}^2({\mathbb{R}^{d}};{\mathbb{R}^{1}}))\bigotimes
M^{2,0}\\([0,T];L_{\rho}^2({\mathbb{R}^{d}};{\mathbb{R}^{d}}))$ and
for an arbitrary $\varphi\in C_c^{0}(\mathbb{R}^d;\mathbb{R}^1)$,
\begin{eqnarray}\label{zhang662}
\int_{\mathbb{R}^{d}}Y_s^{t,x,n}\varphi(x)dx&=&\int_{\mathbb{R}^{d}}h(X_{T}^{t,x})\varphi(x)dx\nonumber\\
&&+\int_{s}^{T}\int_{\mathbb{R}^{d}}f(r,X_{r}^{t,x},Y_r^{t,x,n},Z_r^{t,x,n})\varphi(x)dxdr\nonumber\\
&&-\sum_{j=1}^{n}\int_{s}^{T}\int_{\mathbb{R}^{d}}g_j(r,X_{r}^{t,x},Y_r^{t,x,n},Z_r^{t,x,n})\varphi(x)dxd^\dagger{\hat{\beta}}_j(r)\nonumber\\
&&-\int_{s}^{T}\langle\int_{\mathbb{R}^{d}}Z_r^{t,x,n}\varphi(x)dx,dW_r\rangle\
\ \ P-{\rm a.s.}
\end{eqnarray}
We claim $({Y}_\cdot^{t,\cdot,n},{Z}_\cdot^{t,\cdot,n})$ is a Cauchy
sequence in
$S^{2,0}([0,T];L_{\rho}^2({\mathbb{R}^{d}};{\mathbb{R}^{1}}))\bigotimes
M^{2,0}\\([0,T];L_{\rho}^2({\mathbb{R}^{d}};{\mathbb{R}^{d}}))$. For
this, applying It$\hat {\rm o}$'s formula to ${\rm
e}^{Kr}{{|{Y}_r^{t,x,m}-{Y}_r^{t,x,n}|}^2}$ for a.e.
$x\in\mathbb{R}^{d}$, we have
\begin{eqnarray}\label{zhang667}
&&\int_{\mathbb{R}^{d}}{\rm
e}^{Ks}{|{Y}_s^{t,x,m}-{Y}_s^{t,x,n}|}^2\rho^{-1}(x)dx\nonumber\\
&&+\int_{s}^{T}\int_{\mathbb{R}^{d}}{\rm
e}^{Kr}{|{Y}_r^{t,x,m}-{Y}_r^{t,x,n}|}^2\rho^{-1}(x)dxdr\nonumber\\
&&+\int_{s}^{T}\int_{\mathbb{R}^{d}}{\rm
e}^{Kr}|{Z}_r^{t,x,m}-{Z}_r^{t,x,n}|^2\rho^{-1}(x)dxdr\nonumber\\
&\leq&C_p\sum_{j=n+1}^{m}\{(C_j+\alpha_j)\big(\int_{s}^{T}\int_{\mathbb{R}^{d}}(|{Y}_{r}^{t,x,m}|^2+|{Z}_{r}^{t,x,m}|^2)\rho^{-1}(x)dxdr\nonumber\\
&&\ \ \ \ \ \ \ \ \ \ \ \ \ \ \ \ \ \ \ \ \ \ \ \ \ \ \ +\int_{s}^{T}\int_{\mathbb{R}^{d}}|{g}_j(r,X_{r}^{t,x},0,0)|^2\rho^{-1}(x)dxdr\big)\}\nonumber\\
&&-\sum_{j=1}^{n}\int_{s}^{T}\int_{\mathbb{R}^{d}}2{\rm e}^{Kr}\bar{Y}_r^{t,x,m,n}\bar{g}^{m,n}_j(r,x)\rho^{-1}(x)dxd^\dagger{\hat{\beta}}_j(r)\nonumber\\
&&-\sum_{j=n+1}^{m}\int_{s}^{T}\int_{\mathbb{R}^{d}}2{\rm e}^{Kr}\bar{Y}_r^{t,x,m,n}{g}_j(r,X_{r}^{t,x},{Y}_{r}^{t,x,m},{Z}_{r}^{t,x,m})\rho^{-1}(x)dxd^\dagger{\hat{\beta}}_j(r)\nonumber\\
&&-\int_{s}^{T}\langle\int_{\mathbb{R}^{d}}2{\rm
e}^{Kr}\bar{Y}_r^{t,x,m,n}\bar{Z}_r^{t,x,m,n}\rho^{-1}(x)dx,dW_r\rangle.
\end{eqnarray}
The claim is true by taking expectation and applying Lemma
\ref{qi045} and Proposition \ref{zz58}, as $n$,
$m\longrightarrow\infty$
\begin{eqnarray}\label{zz59}
&&E[\int_{0}^{T}\int_{\mathbb{R}^{d}}{\rm
e}^{Kr}{|{Y}_r^{t,x,m}-{Y}_r^{t,x,n}|}^2\rho^{-1}(x)dxdr]\nonumber\\
&&+E[\int_{0}^{T}\int_{\mathbb{R}^{d}}{\rm
e}^{Kr}|{Z}_r^{t,x,m}-{Z}_r^{t,x,n}|^2\rho^{-1}(x)dxdr]
\longrightarrow0
\end{eqnarray}
and by the B-D-G inequality
\begin{eqnarray*}
E[\sup_{0\leq s\leq T}\int_{\mathbb{R}^{d}}{\rm
e}^{Ks}{|{Y}_r^{t,x,m}-{Y}_r^{t,x,n}|}^2\rho^{-1}(x)dx]
\longrightarrow0.
\end{eqnarray*}
Denote its limit by $({Y}_s^{t,x},{Z}_s^{t,x})$.

We will show that $({Y}_s^{t,x},{Z}_s^{t,x})$ satisfies (\ref{qi22})
for an arbitrary $\varphi\in C_c^{0}(\mathbb{R}^d;\mathbb{R}^1)$.
For this, we prove that along a subsequence (\ref{zhang662}), the
spatial integral form of Eq.(\ref{zhang66100}), converges to
Eq.(\ref{qi22}) in $L^2(\Omega)$ term by term as
$n\longrightarrow\infty$. Here we only show that along a subsequence
\begin{eqnarray*}
E[\
|\int_{s}^{T}\int_{\mathbb{R}^{d}}\big(f(r,X_{r}^{t,x},Y_r^{t,x,n},Z_r^{t,x,n})-f(r,X_{r}^{t,x},Y_r^{t,x},Z_r^{t,x})\big)\varphi(x)dxdr|^2]\longrightarrow0
\end{eqnarray*}
as $n\longrightarrow\infty$. Other items are under the same
conditions as in Section 3 in \cite{zh-zh}, therefore the
convergence can be dealt with similarly. Notice
\begin{eqnarray*}\label{zz17}
&&E[\ |\int_{s}^{T}\int_{\mathbb{R}^{d}}\big(f(r,X_{r}^{t,x},Y_r^{t,x,n},Z_r^{t,x,n})-f(r,X_{r}^{t,x},Y_r^{t,x},Z_r^{t,x})\big)\varphi(x)dxdr|^2]\nonumber\\
&\leq&TE[\int_{s}^{T}\int_{\mathbb{R}^{d}}|f(r,X_{r}^{t,x},Y_r^{t,x,n},Z_r^{t,x,n})-f(r,X_{r}^{t,x},Y_r^{t,x},Z_r^{t,x})|^2\rho^{-1}(x)dxdr\nonumber\\
&&\ \ \ \ \ \ \times\int_{\mathbb{R}^{d}}|\varphi(x)|^2\rho(x)dx]\nonumber\\
&\leq&C_pE[\int_{s}^{T}\int_{\mathbb{R}^{d}}|f(r,X_{r}^{t,x},Y_r^{t,x,n},Z_r^{t,x,n})-f(r,X_{r}^{t,x},Y_r^{t,x},Z_r^{t,x})|^2\rho^{-1}(x)dxdr]\nonumber\\
&\leq&C_pE[\int_{s}^{T}\int_{\mathbb{R}^{d}}|f(r,X_{r}^{t,x},Y_r^{t,x,n},Z_r^{t,x,n})-f(r,X_{r}^{t,x},Y_r^{t,x,n},Z_r^{t,x})|^2\rho^{-1}(x)dxdr]\nonumber\\
&&+C_pE[\int_{s}^{T}\int_{\mathbb{R}^{d}}|f(r,X_{r}^{t,x},Y_r^{t,x,n},Z_r^{t,x})-f(r,X_{r}^{t,x},Y_r^{t,x},Z_r^{t,x})|^2\rho^{-1}(x)dxdr]\nonumber\\
&\leq&C_pE[\int_{s}^{T}\int_{\mathbb{R}^{d}}|Z_r^{t,x,n}-Z_r^{t,x}|^2\rho^{-1}(x)dxdr]\\
&&+C_pE[\int_{s}^{T}\int_{\mathbb{R}^{d}}|f(r,X_{r}^{t,x},Y_r^{t,x,n},Z_r^{t,x})-f(r,X_{r}^{t,x},Y_r^{t,x},Z_r^{t,x})|^2\rho^{-1}(x)dxdr].\nonumber
\end{eqnarray*}
We only need to prove that along a subsequence
\begin{eqnarray}\label{zhang666}
&&E[\int_{s}^{T}\int_{\mathbb{R}^{d}}|f(r,X_{r}^{t,x},Y_r^{t,x,n},Z_r^{t,x})\\
&&\ \ \ \ \ \ \ \ \ \ \ \ \ \ \
-f(r,X_{r}^{t,x},Y_r^{t,x},Z_r^{t,x})|^2\rho^{-1}(x)dxdr]\longrightarrow0\
\ {\rm as}\ n\longrightarrow\infty.\nonumber
\end{eqnarray}
First we will find a subsequence of $\{Y_r^{t,x,n}\}_{n=1}^\infty$
still denoted by $\{Y_r^{t,x,n}\}_{n=1}^\infty$ s.t.
$Y_r^{t,x,n}\longrightarrow Y_r^{t,x}$ for a.e. $r\in[0,T]$,
$x\in\mathbb{R}^{d}$, a.s. $\omega$ and
$E[\int_{0}^{T}\int_{\mathbb{R}^{d}}\sup_n|Y_r^{t,x,n}|^2\rho^{-1}(x)\\dxdr]<\infty$.
For this, from (\ref{zz59}), we know that
$E[\int_{0}^{T}\int_{\mathbb{R}^{d}}|Y_r^{t,x,n}-Y_r^{t,x}|^2\rho^{-1}(x)\\dxdr]\longrightarrow0$.
Therefore we may assume without losing any generality that
$Y_r^{t,x,n}\longrightarrow Y_r^{t,x}$ for a.e. $r\in[0,T]$,
$x\in\mathbb{R}^{d}$, a.s. $\omega$ and extract a subsequence of
$\{Y_r^{t,x,n}\}_{n=1}^\infty$, still denoted by
$\{Y_r^{t,x,n}\}_{n=1}^\infty$, s.t.
\begin{eqnarray*}
\sqrt{E[\int_{0}^{T}\int_{\mathbb{R}^{d}}|Y_r^{t,x,n+1}-Y_r^{t,x,n}|^2\rho^{-1}(x)dxdr]}\leq{1\over{2^n}}.
\end{eqnarray*}
For any $n$,
\begin{eqnarray*}
|Y_r^{t,x,n}|\leq|Y_r^{t,x,1}|+\sum_{i=1}^{n-1}|Y_r^{t,x,i+1}-Y_r^{t,x,i}|\leq|Y_r^{t,x,1}|+\sum_{i=1}^\infty|Y_r^{t,x,i+1}-Y_r^{t,x,i}|.
\end{eqnarray*}
Then by the triangle inequality of the norm, we have
\begin{eqnarray*}
&&\sqrt{E[\int_{0}^{T}\int_{\mathbb{R}^{d}}\sup_n|Y_r^{t,x,n}|^2\rho^{-1}(x)dxdr]}\\
&\leq&\sqrt{E[\int_{0}^{T}\int_{\mathbb{R}^{d}}(|Y_r^{t,x,1}|+\sum_{i=1}^\infty|Y_r^{t,x,i+1}-Y_r^{t,x,i}|)^2\rho^{-1}(x)dxdr]}\\
&\leq&\sqrt{E[\int_{0}^{T}\int_{\mathbb{R}^{d}}|Y_r^{t,x,1}|^2\rho^{-1}(x)dxdr]}\nonumber\\
&&+\sum_{i=1}^\infty\sqrt{E[\int_{0}^{T}\int_{\mathbb{R}^{d}}|Y_r^{t,x,i+1}-Y_r^{t,x,i}|^2\rho^{-1}(x)dxdr]}\\
&\leq&\sqrt{E[\int_{0}^{T}\int_{\mathbb{R}^{d}}|Y_r^{t,x,1}|^2\rho^{-1}(x)dxdr]}+\sum_{i=1}^\infty{1\over{2^i}}\\
&<&\infty.
\end{eqnarray*}
It therefore follows from Condition (H.4) that, for this subsequence
$\{Y_r^{t,x,n}\}_{n=1}^\infty$,
\begin{eqnarray*}
&&E[\int_{0}^{T}\int_{\mathbb{R}^{d}}\sup_n|f(r,X_{r}^{t,x},Y_r^{t,x,n},Z_r^{t,x})-f(r,X_{r}^{t,x},Y_r^{t,x},Z_r^{t,x})|^2\rho^{-1}(x)dxdr]\\
&\leq&C_pE[\int_{0}^{T}\int_{\mathbb{R}^{d}}(1+\sup_n|Y_r^{t,x,n}|^2+|Y_r^{t,x}|^2+|Z_r^{t,x}|^2)\rho^{-1}(x)dxdr]<\infty.
\end{eqnarray*}
Then, (\ref{zhang666}) follows from applying Lebesgue's dominated
convergence theorem and Condition (H.6).
The proof of Theorem \ref{zz48} is completed.
$\hfill\diamond$\\

\subsection{The corresponding SPDEs}\label{s37}
We first consider the following SPDE with finite dimensional noise:
\begin{eqnarray}\label{zz12}
u^n(t,x)&=&h(x)+\int_{t}^{T}[\mathscr{L}u^n(s,x)+f\big(s,x,u^n(s,x),(\sigma^*\nabla u^n)(s,x)\big)]ds\\
&&-\sum_{j=1}^{n}\int_{t}^{T}g_j\big(s,x,u^n(s,x),(\sigma^*\nabla
u^n)(s,x)\big)d^\dagger{\hat{\beta}}_j(s),\ \ \ \ 0\leq t\leq s\leq
T.\nonumber
\end{eqnarray}
In the previous subsection, we proved the existence and uniqueness
of solution of BDSDE (\ref{qi20}) and obtained the solution
$({Y}_s^{t,x},{Z}_s^{t,x})$ by taking the limit of
$({Y}_s^{t,x,n},{Z}_s^{t,x,n})$ of the solutions of
Eq.(\ref{zhang66100}) in the space
$S^{2,0}([0,T];L_{\rho}^2({\mathbb{R}^{d}};{\mathbb{R}^{1}}))\bigotimes\\
M^{2,0}([0,T];L_{\rho}^2({\mathbb{R}^{d}};{\mathbb{R}^{d}}))$ along
a subsequence. We still start from Eq.(\ref{zhang66100}) in this
subsection.
\begin{prop}\label{zz004}
Under Conditions {\rm(H.1)}--{\rm(H.7)}, assume
Eq.(\ref{zhang66100}) has a unique solution
$(Y_r^{t,x,n},\\Z_r^{t,x,n})$, then for any $t\leq s\leq T$,
\begin{eqnarray*}
Y_r^{s,X_s^{t,x},n}=Y_r^{t,x,n}\ {\rm and}\
Z_r^{s,X_s^{t,x},n}=Z_r^{t,x,n}\ {\rm for}\ r\in[s,T],\ {\rm a.a.}\
x\in{\mathbb{R}^{d}}\ {\rm a.s.}
\end{eqnarray*}
\end{prop}
{\em Proof}. The proof is similar to the proof of Proposition 3.4 in
\cite{zh-zh}. Here Lemma \ref{zz002} plays the same role as Lemma
3.3 in that proof.
$\hfill\diamond$\\

A direct application of Proposition \ref{zz004} and Fubini theorem
immediately leads to
\begin{prop}\label{zz005} Under Conditions {\rm(H.1)}--{\rm(H.7)}, if we define $u^n(t,x)=Y_t^{t,x,n}$, $v^n(t,x)=Z_t^{t,x,n}$, then
\begin{eqnarray*}
u^n(s,X^{t,x}_s)=Y_s^{t,x,n},\ v^n(s,X^{t,x}_s)=Z_s^{t,x,n}\ {\rm
for}\ {\rm a.a.}\ s\in[t,T],\ x\in{\mathbb{R}^{d}}\ {\rm a.s.}
\end{eqnarray*}
\end{prop}
\begin{thm}\label{zz006} Under Conditions {\rm(H.1)}--{\rm(H.7)},
if we define $u^n(t,x)=Y_t^{t,x,n}$, where
$(Y_s^{t,x,n},Z_s^{t,x,n})$ is the solution of
Eq.(\ref{zhang66100}), then $u^n(t,x)$ is the unique weak solution
of Eq.(\ref{zz12}). Moreover,
\begin{eqnarray*}
u^n(s,X_s^{t,x})=Y_s^{t,x,n},\ (\sigma^*\nabla
u^n)(s,X^{t,x}_s)=Z_s^{t,x,n}\ {\rm for}\ {\rm a.a.}\ s\in[t,T],\
x\in\mathbb{R}^{d}\ {\rm a.s.}
\end{eqnarray*}
\end{thm}
{\em Proof}. \underline{Uniqueness}. Let $u^n$ be a solution of
Eq.(\ref{zz12}). Define
\begin{eqnarray*}
&&F^n(s,x)=f\big(s,x,u^n(s,x),(\sigma^*\nabla u^n)(s,x)\big),\\
&&G^n_j(s,x)=g_j\big(s,x,u^n(s,x),(\sigma^*\nabla u^n)(s,x)\big).
\end{eqnarray*}
Since $u^n$ is the solution, so
$E[\int_{0}^{T}\int_{\mathbb{R}^d}\big(|u^n(s,x)|^2+|(\sigma^*\nabla
u^n)(s,x)|^2\big)\rho^{-1}(x)dxds]\\<\infty$ and
\begin{eqnarray}\label{zz14}
&&E[\int_{0}^{T}\int_{\mathbb{R}^d}(|F^n(s,x)|^2+\sum_{j=1}^n|G^n_j(s,x)|^2)\rho^{-1}(x)dxds]\nonumber\\
&=&E[\int_{0}^{T}\int_{\mathbb{R}^d}\big(|f\big(s,x,u^n(s,x),(\sigma^*\nabla u^n)(s,x)\big)|^2\nonumber\\
&&\ \ \ \ \ \ \ \ \ \ \ \ \ \ \ +\sum_{j=1}^n|g_j\big(s,x,u^n(s,x),(\sigma^*\nabla u^n)(s,x)\big)|^2\big)\rho^{-1}(x)dxds]\nonumber\\
&\leq&C_pE[\int_{0}^{T}\int_{\mathbb{R}^d}\big(1+|u^n(s,x)|^2+|(\sigma^*\nabla u^n)(s,x)|^2\nonumber\\
&&\ \ \ \ \ \ \ \ \ \ \ \ \ \ \ \ \ \ \ +\sum_{j=1}^n|g_j(s,x,0,0)|^2\big)\rho^{-1}(x)dxds]\nonumber\\
&<&\infty.
\end{eqnarray}
If we define $Y_s^{t,x,n}=u^n(s,X_s^{t,x})$ and
$Z_s^{t,x,n}=(\sigma^*\nabla u^n)(s,X^{t,x}_s)$, then by Lemma
\ref{qi045},
\begin{eqnarray*}
&&E[\int_{t}^{T}\int_{\mathbb{R}^d}(|Y_s^{t,x,n}|^2+|Z_s^{t,x,n}|^2)\rho^{-1}(x)dxds]\nonumber\\
&\leq&C_pE[\int_{t}^{T}\int_{\mathbb{R}^d}|u^n(s,x)|^2+|(\sigma^*\nabla
u^n)(s,x)|^2\rho^{-1}(x)dxds]<\infty.
\end{eqnarray*}
Using some ideas of Theorem 2.1 in \cite{ba-ma}, similar to the
argument as in Section 4 in \cite{zh-zh}, we have for $t\leq s\leq
T$, $(Y_\cdot^{t,\cdot,n},Z_\cdot^{t,\cdot,n})\in M^{2,0}\\([t,T];
L_{\rho}^2({\mathbb{R}^{d}};{\mathbb{R}^{1}}))\bigotimes
M^{2,0}([t,T];L_{\rho}^2({\mathbb{R}^{d}};{\mathbb{R}^{d}}))$ solves
the following BDSDE:
\begin{eqnarray}\label{zz15}
Y_s^{t,x,n}&=&h(X_{T}^{t,x})+\int_{s}^{T}F^n(r,X_{r}^{t,x})dr\nonumber\\
&&-\sum_{j=1}^{n}\int_{s}^{T}G^n_j(r,X_{r}^{t,x})d^\dagger{\hat{\beta}}_j(r)-\int_{s}^{T}\langle
Z^{t,x,n}_r,dW_r\rangle.
\end{eqnarray}
Multiply $\varphi\in C_c^{0}(\mathbb{R}^d;\mathbb{R}^1)$ on both
sides and then take the integration over $\mathbb{R}^d$. Noting the
definition of $F^n(s,x)$, $G^n_j(s,x)$, $Y_s^{t,x,n}$ and
$Z_s^{t,x,n}$, we have that $(Y_s^{t,x,n},Z_s^{t,x,n})$ satisfies
the spatial integration form of Eq.(\ref{zhang66100}). By Corollary
\ref{zz003}, $Y_\cdot^{t,\cdot,n}\in
S^{2,0}([t,T];L_{\rho}^2({\mathbb{R}^{d}};{\mathbb{R}^{1}}))$ and
therefore $(Y_s^{t,x,n},Z_s^{t,x,n})$ is a solution of
Eq.(\ref{zhang66100}). If there is another solution $\hat{u}$ to
Eq.(\ref{zz12}), then by the same procedure, we can find another
solution $(\hat{Y}_s^{t,x,n},\hat{Z}_s^{t,x,n})$ to
Eq.(\ref{zhang66100}), where
\begin{eqnarray*}
\hat{Y}_s^{t,x,n}=\hat{u}^n(s,X_s^{t,x})\ {\rm and}\
\hat{Z}_s^{t,x,n}=(\sigma^*\nabla \hat{u}^n)(s,X^{t,x}_s).
\end{eqnarray*}
By Theorem \ref{zz000}, the solution of Eq.(\ref{zhang66100}) is
unique, therefore
\begin{eqnarray*}
Y_s^{t,x,n}=\hat{Y}_s^{t,x,n}\ {\rm for}\ {\rm a.a.}\ s\in[t,T],\
x\in\mathbb{R}^{d}\ {\rm a.s.}
\end{eqnarray*}
Especially for $t=0$,
\begin{eqnarray*}
Y_s^{0,x,n}=\hat{Y}_s^{0,x,n}\ {\rm for}\ {\rm a.a.}\ s\in[0,T],\
x\in\mathbb{R}^{d}\ {\rm a.s.}
\end{eqnarray*}
By Lemma \ref{qi045} again,
\begin{eqnarray*}
&&E[\int_{0}^{T}\int_{\mathbb{R}^d}|u^n(s,x)-\hat{u}^n(s,x)|^2\rho^{-1}(x)dxds]\nonumber\\
&\leq&C_pE[\int_{0}^{T}\int_{\mathbb{R}^d}|u^n(s,X^{0,x}_s)-\hat{u}^n(s,X^{0,x}_s)|^2\rho^{-1}(x)dxds]\nonumber\\
&=&C_pE[\int_{0}^{T}\int_{\mathbb{R}^d}|Y_s^{0,x,n}-\hat{Y}_s^{0,x,n}|^2)\rho^{-1}(x)dxds]\nonumber\\
&=&0.
\end{eqnarray*}
So $u^n(s,x)=\hat{u}^n(s,x)$ for a.a. $s\in[0,T]$,
$x\in\mathbb{R}^{d}$ a.s. The uniqueness is proved.
\\

\underline{Existence}. For each $(t,x)\in[0,T]\otimes\mathbb{R}^d$,
define $u^n(t,x)=Y_t^{t,x,n}$ and $v^n(t,x)=Z_t^{t,x,n}$, where
$({Y}_\cdot^{t,\cdot,n},{Z}_\cdot^{t,\cdot,n})\in
S^{2,0}([0,T];L_{\rho}^2({\mathbb{R}^{d}};{\mathbb{R}^{1}}))\bigotimes
M^{2,0}([0,T];\\L_{\rho}^2({\mathbb{R}^{d}};{\mathbb{R}^{d}}))$ is
the solution of Eq.(\ref{zhang66100}). Then by Proposition
\ref{zz005},
\begin{eqnarray*}
u^n(s,X^{t,x}_s)=Y_s^{t,x,n},\ v^n(s,X^{t,x}_s)=Z_s^{t,x,n}\ {\rm
for}\ {\rm a.a.}\ s\in[t,T],\ x\in{\mathbb{R}^{d}}\ {\rm a.s.}
\end{eqnarray*}
Set
\begin{eqnarray*}
&&F^n(s,x)=f\big(s,x,u^n(s,x),v^n(s,x)\big),\\
&&G^n_j(s,x)=g_j\big(s,x,u^n(s,x),v^n(s,x)\big).
\end{eqnarray*}
Then it is easy to see that
$(Y_\cdot^{t,\cdot,n},Z_\cdot^{t,\cdot,n})\in M^{2,0}([t,T];
L_{\rho}^2({\mathbb{R}^{d}};{\mathbb{R}^{1}}))\bigotimes
M^{2,0}\\([t,T];L_{\rho}^2({\mathbb{R}^{d}};{\mathbb{R}^{d}}))$ is a
solution of Eq.(\ref{zz15}) with above $F^n$ and $G_j^N$. Moreover,
by Lemma \ref{qi045},
\begin{eqnarray*}
E[\int_{0}^{T}\int_{\mathbb{R}^d}|u^n(s,x)|^2+|v^n(s,x)|^2\rho^{-1}(x)dxds]<\infty.
\end{eqnarray*}
Then from a similar computation as in (\ref{zz14}) we have
\begin{eqnarray*}
E[\int_{0}^{T}\int_{\mathbb{R}^d}(|F^n(s,x)|^2+\sum_{j=1}^n|G^n_j(s,x)|^2)\rho^{-1}(x)dxds]<\infty.
\end{eqnarray*}
Now using some ideas of Theorem 2.1 in \cite{ba-ma}, similar to the
argument as in Section 4 in \cite{zh-zh}, we know that
$v^n(s,x)=(\sigma^*\nabla u^n)(s,x)$ and $u^n$ is the weak solution
of the following SPDE:
\begin{eqnarray}\label{zz16}
u^n(t,x)&=&h(x)+\int_{t}^{T}[\mathscr{L}u^n(s,x)+F^n(s,x)]ds\nonumber\\
&&-\sum_{j=1}^{n}\int_{t}^{T}G^n_j(s,x)d^\dagger{\hat{\beta}}_j(s),\
\ \ \ 0\leq t\leq s\leq T.
\end{eqnarray}
Noting the definition of $F^n(s,x)$ and $G^n_j(s,x)$ and the fact
that $v^n(s,x)=(\sigma^*\nabla u^n)(s,x)$, from (\ref{zz16}), we
have that $u^n$
is the weak solution of Eq.(\ref{zz12}). $\hfill\diamond$\\

In the rest part of this subsection, we study Eq.(\ref{zhang685})
with $f$ and $g$ allowed to depend on time.
If $(Y_s^{t,x},Z_s^{t,x})$ is the solution of Eq.(\ref{qi20}) and we
define $u(t,x)=Y_t^{t,x}$, then by Proposition 4.2 in \cite{zh-zh},
we have $\sigma^*\nabla u(t,x)$ exists for a.a. $t\in[0,T]$,
$x\in\mathbb{R}^{d}$ a.s., and
\begin{eqnarray}\label{zhang660}
u(s,X_s^{t,x})=Y_s^{t,x},\ (\sigma^*\nabla
u)(s,X^{t,x}_s)=Z_s^{t,x}\ {\rm for}\ {\rm a.a.}\ s\in[t,T],\
x\in\mathbb{R}^{d}\ {\rm a.s.}\nonumber\\
\end{eqnarray}
Also by Theorem \ref{zz006} and Lemma \ref{qi045}, we have
\begin{eqnarray}\label{zz18}
&&E[\int_{0}^{T}\int_{\mathbb{R}^d}|u^n(s,x)-u(s,x)|^2\rho^{-1}(x)dxds]\nonumber\\
&&+E[\int_{0}^{T}\int_{\mathbb{R}^d}|(\sigma^*\nabla u^n)(s,x)-(\sigma^*\nabla u)(s,x)|^2\rho^{-1}(x)dxds]\nonumber\\
&\leq&C_pE[\int_{0}^{T}\int_{\mathbb{R}^d}|u^n(s,X^{0,x}_s)-u(s,X^{0,x}_s)|^2\rho^{-1}(x)dxds]\nonumber\\
&&+C_pE[\int_{0}^{T}\int_{\mathbb{R}^d}|(\sigma^*\nabla u^n)(s,X^{0,x}_s)-(\sigma^*\nabla u)(s,X^{0,x}_s)|^2\rho^{-1}(x)dxds]\nonumber\\
&=&C_pE[\int_{0}^{T}\int_{\mathbb{R}^d}|{Y}_s^{0,x,n}-{Y}_s^{0,x}|^2\rho^{-1}(x)dxds]\nonumber\\
&&+C_pE[\int_{0}^{T}\int_{\mathbb{R}^d}|{Z}_s^{0,x,n}-{Z}_s^{0,x}|^2\rho^{-1}(x)dxds]\longrightarrow0,\
\ {\rm as}\ n\rightarrow\infty.
\end{eqnarray}
With (\ref{zz18}), we prove the other main theorem in this section.
\\

{\em Proof of Theorem \ref{zz45}}. We only need to verify that this
$u$ defined through $Y_t^{t,x}$ is the unique weak solution of
Eq.(\ref{zhang685}). By Lemma \ref{qi045} and (\ref{zhang660}), it
is easy to see that
\begin{eqnarray*}
(\sigma^*\nabla u)(t,x)=Z_t^{t,x}\ {\rm for}\ {\rm a.a.}\
t\in[0,T],\ x\in\mathbb{R}^d\ {\rm a.s.}
\end{eqnarray*}
Furthermore, using the generalized equivalence norm principle again
we have
\begin{eqnarray}\label{zz46}
&&E[\int_{0}^{T}\int_{\mathbb{R}^d}(|u(s,x)|^2+|(\sigma^*\nabla u)(s,x)|^2)\rho^{-1}(x)dxds]\nonumber\\
&\leq&C_pE[\int_{0}^{T}\int_{\mathbb{R}^d}(|u(s,X^{0,x}_s)|^2+|(\sigma^*\nabla u)(s,X^{0,x}_s)|^2)\rho^{-1}(x)dxds]\nonumber\\
&=&C_pE[\int_{0}^{T}\int_{\mathbb{R}^d}(|Y_s^{0,x}|^2+|Z_s^{0,x}|^2)\rho^{-1}(x)dxds]<\infty.
\end{eqnarray}
Now we will verify that $u(t,x)$ satisfies (\ref{qi16}). Since
$u^n(t,x)$ is the weak solution of SPDE (\ref{zz12}), so for any
$\varphi\in C_c^{\infty}(\mathbb{R}^d;\mathbb{R}^1)$, $u^n(t,x)$
satisfies
\begin{eqnarray}\label{qi25}
&&\int_{\mathbb{R}^{d}}u^n(t,x)\varphi(x)dx-\int_{\mathbb{R}^{d}}h(x)\varphi(x)dx\nonumber\\
&&-{1\over2}\int_{t}^{T}\int_{\mathbb{R}^{d}}(\sigma^*\nabla u^n)(s,x)(\sigma^*\nabla\varphi)(x)dxds\nonumber\\
&&-\int_{t}^{T}\int_{\mathbb{R}^{d}}u^n(s,x)\nabla\big((b-\tilde{A})\varphi\big)(x)dxds\nonumber\\
&=&\int_{t}^{T}\int_{\mathbb{R}^{d}}f\big(s,x,u^n(s,x),(\sigma^*\nabla u^n)(s,x)\big)\varphi(x)dxds\\
&&-\sum_{j=1}^{n}\int_{t}^{T}\int_{\mathbb{R}^{d}}g_j\big(s,x,u^n(s,x),(\sigma^*\nabla
u^n)(s,x)\big)\varphi(x)dxd^\dagger{\hat{\beta}}_j(s)\ \ \ \ P-{\rm
a.s.}\nonumber
\end{eqnarray}
By proving that along a subsequence (\ref{qi25}) converges to
(\ref{qi16}) in $L^2(\Omega)$, we have that $u(t,x)$ satisfies
(\ref{qi16}). We only need to show that along a sequence as
$n\longrightarrow\infty$,
\begin{eqnarray*}
&&E[\
|\int_{t}^{T}\int_{\mathbb{R}^{d}}\big(f(s,x,u^n(s,x),(\sigma^*\nabla
u^n)(s,x))\nonumber\\
&&\ \ \ \ \ \ \ \ \ \ \ \ \ \ \ \ \ \ -f(s,x,u(s,x),(\sigma^*\nabla
u)(s,x))\big)\varphi(x)dxds|^2]\longrightarrow0.
\end{eqnarray*}
First note
\begin{eqnarray*}
&&E[\ |\int_{t}^{T}\int_{\mathbb{R}^{d}}\big(f(s,x,u^n(s,x),(\sigma^*\nabla u^n)(s,x))\nonumber\\
&&\ \ \ \ \ \ \ \ \ \ \ \ \ \ \ \ \ \ -f(s,x,u(s,x),(\sigma^*\nabla u)(s,x))\big)\varphi(x)dxds|^2]\\
&\leq&C_pE[\int_{t}^{T}\int_{\mathbb{R}^{d}}|f(s,x,u^n(s,x),(\sigma^*\nabla u^n)(s,x))\nonumber\\
&&\ \ \ \ \ \ \ \ \ \ \ \ \ \ \ \ \ \ -f(s,x,u(s,x),(\sigma^*\nabla u)(s,x))|^2\rho^{-1}(x)dxds]\\
&\leq&C_pE[\int_{t}^{T}\int_{\mathbb{R}^{d}}|f(s,x,u^n(s,x),(\sigma^*\nabla u^n)(s,x))\nonumber\\
&&\ \ \ \ \ \ \ \ \ \ \ \ \ \ \ \ \ \ -f(s,x,u^n(s,x),(\sigma^*\nabla u)(s,x))|^2\rho^{-1}(x)dxds]\\
&&+C_pE[\int_{t}^{T}\int_{\mathbb{R}^{d}}|f(s,x,u^n(s,x),(\sigma^*\nabla u)(s,x))\nonumber\\
&&\ \ \ \ \ \ \ \ \ \ \ \ \ \ \ \ \ \ \ \ -f(s,x,u(s,x),(\sigma^*\nabla u)(s,x))|^2\rho^{-1}(x)dxds]\\
&\leq&C_pE[\int_{t}^{T}\int_{\mathbb{R}^{d}}|(\sigma^*\nabla u^n)(s,x)-(\sigma^*\nabla u)(s,x)|^2\rho^{-1}(x)dxds]\\
&&+C_pE[\int_{t}^{T}\int_{\mathbb{R}^{d}}|f(s,x,u^n(s,x),(\sigma^*\nabla
u)(s,x))\nonumber\\
&&\ \ \ \ \ \ \ \ \ \ \ \ \ \ \ \ \ \ \ \
-f(s,x,u(s,x),(\sigma^*\nabla u)(s,x))|^2\rho^{-1}(x)dxds].
\end{eqnarray*}
We face a similar situation as in (\ref{zz17}) and only need to
prove that along a subsequence as $n\longrightarrow\infty$,
\begin{eqnarray}\label{zhang664}
&&E[\int_{t}^{T}\int_{\mathbb{R}^{d}}|f(s,x,u^n(s,x),(\sigma^*\nabla
u)(s,x))\\
&&\ \ \ \ \ \ \ \ \ \ \ \ \ \ \ \ \ \ -f(s,x,u(s,x),(\sigma^*\nabla
u)(s,x))|^2\rho^{-1}(x)dxds]\longrightarrow0.\nonumber
\end{eqnarray}
For this, note that we have (\ref{zz18}) which plays the same role
as ({\ref{zz59}) in the proof of Theorem \ref{zz48}. Thus we can
find a subsequence of $\{u^n(s,x)\}_{n=1}^\infty$ still denoted by
$\{u^n(s,x)\}_{n=1}^\infty$ s.t. $u^n(s,x)\longrightarrow u(s,x)$
for a.e. $s\in[0,T]$, $x\in\mathbb{R}^{d}$, a.s. $\omega$ and
$E[\int_{0}^{T}\int_{\mathbb{R}^{d}}\sup_n|u^n(s,x)|^2\rho^{-1}(x)dxds]<\infty$.
It turns out that, for this subsequence $\{u^n(s,x)\}_{n=1}^\infty$,
by Condition (H.4), we have
\begin{eqnarray*}
&&E[\int_{0}^{T}\int_{\mathbb{R}^{d}}\sup_n|f(s,x,u^n(s,x),(\sigma^*\nabla
u)(s,x))\nonumber\\
&&\ \ \ \ \ \ \ \ \ \ \ \ \ \ \ \ \ \ \
-f(s,x,u(s,x),(\sigma^*\nabla u)(s,x))|^2\rho^{-1}(x)dxds]<\infty.
\end{eqnarray*}
Thus (\ref{zhang664}) follows from using Lebesgue's dominated
convergence theorem. Convergences of other terms in (\ref{qi25}) are
easy to check.

Therefore $u(t,x)$ satisfies (\ref{qi16}), i.e. it is a weak
solution of Eq.(\ref{zhang685}) with $u(T,x)=h(x)$. We can prove the
uniqueness following a similar argument in Theorem \ref{zz006}. $\hfill\diamond$\\

\section{Stationary Solutions of SPDEs and Infinite Horizon BDSDEs}\label{s38}

\setcounter{equation}{0}

In this section, first we will give the proof of Theorem
\ref{qi888}. Then we show that the conditions in Theorem \ref{qi888}
are satisfied, i.e. both Theorem \ref{qi043} and Theorem \ref{qi044}
are true under our assumptions.
\subsection{Proof of Theorem \ref{qi888}}
{\em Proof}. First note that Eq.(\ref{qi13}) is equivalent to the
following BDSDE
\begin{eqnarray}\label{qi8}
\left\{\begin{array}{l} Y_{s}^{t,x}=Y_{T}^{t,x}+\int_{s}^{T}f(X_r^{t,x},Y_r^{t,x}, Z_r^{t,x})dr\\
\ \ \ \ \ \ \ \ \ \ -\int_{s}^{T}g(X_r^{t,x},Y_r^{t,x}, Z_r^{t,x})d^\dagger{\hat{B}}_r -\int_{s}^{T}Z_r^{t,x}dW_r\\
\lim_{T\rightarrow\infty}{\rm e}^{-KT}Y_{T}=0 \ \ \ {\rm a.s.}
\end{array}\right.
\end{eqnarray}
Let ${B}_ u=\hat{B}_{T'-u}-\hat{B}_{T'}$ for arbitrary $T'>0$ and
$-\infty<u\leq T'$. Then ${B}_u$ is a Brownian motion with
${B}_0=0$. For any $r\geq0$, applying $\hat{\theta}_r$ on ${B}_u$,
we have
\begin{eqnarray*}
\hat{\theta}_r\circ{B}_u&=&\hat{\theta}_r\circ(\hat{B}_{T'-u}-\hat{B}_{T'})=\hat{B}_{T'-u+r}-\hat{B}_{T'+r}\\
&=&(\hat{B}_{T'-u+r}-\hat{B}_{T'})-(\hat{B}_{T'+r}-\hat{B}_T')={B}_{u-r}-{B}_{-r}.
\end{eqnarray*}
So for $0\leq s\leq T\leq T'$ and $\{h(u,\cdot)\}_{u\geq0}$ being a
$\mathscr{F}_u$-measurable and locally square integrable stochastic
process with values on
${\mathcal{L}^2_{U_0}(L_{\rho}^2({\mathbb{R}^{d}};{\mathbb{R}^{1}}))}$,
we have the relationship between the forward integral and backward
It$\hat {\rm o}$ integral (c.f. \cite{zh-zh})
\begin{eqnarray*}
\int_s^Th(u,\cdot)d^\dagger
B_u=-\int_{T'-T}^{T'-s}h(T'-u,\cdot)dB_u\ \ \ {\rm a.s.}
\end{eqnarray*}
and 
for arbitrary $T\geq0$, $0\leq s\leq T$, $r\geq0$,
\begin{eqnarray}\label{qi5}
\hat{\theta}_r\circ\int_{s}^{T}h(u,\cdot)d^\dagger
\hat{B}_u=\int_{s+r}^{T+r}\hat{\theta}_r\circ
h(u-r,\cdot)d^\dagger{\hat{B}}_u.
\end{eqnarray}
Therefore for a.e. $x\in\mathbb{R}^{d}$,
\begin{eqnarray*}\label{qi5xx}
\hat{\theta}_r\circ\int_{s}^{T}h(u,x)d^\dagger
\hat{B}_u=\int_{s+r}^{T+r}\hat{\theta}_r\circ
h(u-r,x)d^\dagger{\hat{B}}_u.
\end{eqnarray*}
Since $(Y_{\cdot}^{t,\cdot},Z_{\cdot}^{t,\cdot})\in S^{2,-K}\bigcap
M^{2,-K}([0,\infty);L_{\rho}^2({\mathbb{R}^{d}};{\mathbb{R}^{1}}))
\bigotimes
M^{2,-K}([0,\infty);L_{\rho}^2\\({\mathbb{R}^{d}};{\mathbb{R}^{d}}))$
is the unique solution of Eq.(\ref{qi13}), it follows that
$g(X_\cdot^{t,\cdot},Y_\cdot^{t,\cdot}, Z_\cdot^{t,\cdot})$ is
locally square integrable with values on
${\mathcal{L}^2_{U_0}(L_{\rho}^2({\mathbb{R}^{d}};{\mathbb{R}^{1}}))}$.
Therefore by (\ref{qi18}) and (\ref{qi5}), for a.e.
$x\in\mathbb{R}^{d}$
\begin{eqnarray}\label{qi6}
&&\hat{\theta}_r\circ\int_{s}^{T}g(X_u^{t,x},Y_u^{t,x},
Z_u^{t,x})d^\dagger
\hat{B}_u\nonumber\\
&=&\hat{\theta}_r\circ\int_{s+r}^{T+r}
g(X_{u-r}^{t,x},Y_{u-r}^{t,x},Z_{u-r}^{t,x})d^\dagger{\hat{B}}_u\nonumber\\
&=&\int_{s+r}^{T+r}g(X_{u}^{t+r,x},\hat{\theta}_r\circ
Y_{u-r}^{t,x},\hat{\theta}_r\circ
Z_{u-r}^{t,x})d^\dagger{\hat{B}}_u.
\end{eqnarray}
Now applying the operator $\hat{\theta}_r$ on both sides of
Eq.(\ref{qi8}) and by (\ref{qi6}), we know that $\hat{\theta}_r\circ
Y_s^{t,x}$ satisfies the following equation
\begin{eqnarray}
\left\{\begin{array}{l}\label{qi7} \hat{\theta}_r\circ
Y_s^{t,x}=\hat{\theta}_r\circ
Y_T^{t,x}+\int_{s+r}^{T+r}f(X_{u}^{t+r,x},\hat{\theta}_r\circ
Y_{u-r}^{t,x},\hat{\theta}_r\circ
Z_{u-r}^{t,x})du\\
\ \ \ \ \ \ \ \ \ \ \ \ \ \ \ \
-\int_{s+r}^{T+r}g(X_{u}^{t+r,x},\hat{\theta}_r\circ
Y_{u-r}^{t,x},\hat{\theta}_r\circ
Z_{u-r}^{t,x})d^\dagger{\hat{B}}_u\\
\ \ \ \ \ \ \ \ \ \ \ \ \ \ \ \
-\int_{s+r}^{T+r}\hat{\theta}_r\circ Z_{u-r}^{t,x}dW_u\\
\lim_{T\rightarrow\infty}{\rm e}^{-K(T+r)}(\hat{\theta}_r\circ
Y_T^{t,x})=0\ \ \ {\rm a.s.}
\end{array}\right.
\end{eqnarray}
On the other hand, from Eq.(\ref{qi8}) it is obvious that
\begin{eqnarray}\label{qi2}
\left\{\begin{array}{l}
Y_{s+r}^{t+r,x}=Y_{T+r}^{t+r,x}+\int_{s+r}^{T+r}f(X_{u}^{t+r,x},Y_{u}^{t+r,x},Z_{u}^{t+r,x})du\\
\ \ \ \ \ \ \ \ \ \ \ \ \ -\int_{s+r}^{T+r}g(X_{u}^{t+r,x},Y_{u}^{t+r,x},Z_{u}^{t+r,x})d^\dagger{\hat{B}}_u\\
\ \ \ \ \ \ \ \ \ \ \ \ \ -\int_{s+r}^{T+r}Z_u^{t+r,x}dW_u\\
\lim_{T\rightarrow\infty}{\rm e}^{-K(T+r)}Y_{T+r}^{t+r,x}=0 \ \ \
{\rm a.s.}
\end{array}\right.
\end{eqnarray}
Let $\hat{Y}_{\cdot}^{t,\cdot}=\hat{\theta}_r\circ
Y_{\cdot-r}^{t-r,\cdot}$,
$\hat{Z}_{\cdot}^{t,\cdot}=\hat{\theta}_r\circ
Z_{\cdot-r}^{t-r,\cdot}$. By the uniqueness of the solution of
Eq.(\ref{qi13}) in the space $S^{2,-K}\bigcap
M^{2,-K}([0,\infty);L_{\rho}^2({\mathbb{R}^{d}};{\mathbb{R}^{1}}))
\bigotimes
M^{2,-K}([0,\infty);L_{\rho}^2\\({\mathbb{R}^{d}};{\mathbb{R}^{d}}))$,
it follows from comparing (\ref{qi7}) with (\ref{qi2}) that for any
$r\geq0$ and $t\geq0$, in the space
$L_{\rho}^2({\mathbb{R}^{d}};{\mathbb{R}^{1}})\bigotimes
L_{\rho}^2({\mathbb{R}^{d}};{\mathbb{R}^{d}})$
\begin{eqnarray*}
\hat{\theta}_r\circ
Y_s^{t,\cdot}=\hat{Y}_{s+r}^{t+r,\cdot}=Y_{s+r}^{t+r,\cdot}, \ \
\hat{\theta}_r\circ
Z_s^{t,\cdot}=\hat{Z}_{s+r}^{t+r,\cdot}=Z_{s+r}^{t+r,\cdot}\ \ {\rm
for}\ {\rm all}\ s\geq t\ {\rm a.s.}
\end{eqnarray*}
Then by the perfection procedure (\cite{ar}, \cite{ar-sc}), we can
prove above identities are true for all $s\geq t$, $r\geq0$, but
fixed $t\geq0$ a.s. In particular, for any $t\geq0$, in the space
$L_{\rho}^2({\mathbb{R}^{d}};{\mathbb{R}^{1}})\bigotimes
L_{\rho}^2({\mathbb{R}^{d}};{\mathbb{R}^{d}})$
\begin{eqnarray}\label{qi19}
\hat{\theta}_r\circ Y_{t}^{t,\cdot}=Y_{t+r}^{t+r,\cdot}\ \ \ {\rm
for}\ {\rm all}\ r\geq0\ {\rm a.s.}
\end{eqnarray}
From the assumptions, we also know that $u(t,\cdot)\triangleq
Y_{t}^{t,\cdot}$ is the continuous weak solution of
Eq.(\ref{zhang685}). So we get from (\ref{qi19}) that for any
$t\geq0$, in the space
$L_{\rho}^2({\mathbb{R}^{d}};{\mathbb{R}^{1}})\bigotimes
L_{\rho}^2({\mathbb{R}^{d}};{\mathbb{R}^{d}})$
\begin{eqnarray*}
\hat{\theta}_r\circ u({t,\cdot})=u({t+r,\cdot})\ \ \ {\rm for}\ {\rm
all}\ r\geq0\ {\rm a.s.}
\end{eqnarray*}
Until now, we know "crude" stationary property for $u({t,\cdot})$,
but due to the continuity of $u(t,\cdot)$ w.r.t. $t$ we can obtain
an indistinguishable version of $u({t,\cdot})$, still denoted by
$u({t,\cdot})$, s.t.
\begin{eqnarray*}
\hat{\theta}_{r}\circ u({t,\cdot})=u({t+r,\cdot})\ \ \ {\rm for}\
{\rm all}\ t,\ r\geq0\ \ {\rm a.s.}
\end{eqnarray*}
So we proved that $u(t,\cdot)$ is a ``perfect" stationary weak
solution of Eq.(\ref{zhang685}).

By Definition \ref{qi042}, it follows that
$g\big(\cdot,u(s,\cdot),(\sigma^*\nabla
u)(s,\cdot)\big)\in\mathcal{L}^2_{U_0}(L_{\rho}^2({\mathbb{R}^{d}};{\mathbb{R}^{1}}))$
should be locally square integrable. Now we consider
Eq.(\ref{zhao21}) with cylindrical Brownian motion $B$ on $U_0$. For
arbitrary $T'>0$, let $Y$ be the solution of Eq.(\ref{qi13}) and
$u(t,\cdot)=Y_t^{t,\cdot}$ be the stationary solution of
Eq.(\ref{zhang685}) with $\hat{B}$ chosen as the time reversal of
$B$ from time $T'$, i.e. $\hat{B}_s=B_{T'-s}-B_{T'}$ for $s\geq0$.
Doing the integral transformation in the integration form
(\ref{qi16}) of Eq.(\ref{zhang685}), it is easy to see that
$v_t(x)\triangleq u(T'-t,x)$ satisfies (\ref{zhao21}).

In fact, we can prove a claim that
$v_t(\cdot)(\omega)=Y_{T'-t}^{T'-t,\cdot}(\hat{\omega})$ does not
depend on the choice of $T'$. For this, we only need to show that
for any ${T'}^*\geq T'$,
$Y_{T'-t}^{T'-t,\cdot}(\hat{\omega})=Y_{{T'}^*-t}^{{T'}^*-t,\cdot}({\hat{\omega}}^*)$
when $0\leq t\leq T'$, where $\hat{\omega}(s)={B}_{T'-s}-{B}_{T'}$
and ${\hat{\omega}}^*(s)={B}_{{T'}^*-s}-{B}_{{T'}^*}$. Let
$\hat{\theta}_\cdot$ and $\hat{\theta}^*_\cdot$ be the shifts of
$\hat{\omega}(\cdot)$ and ${\hat{\omega}}^*(\cdot)$ respectively.
Since by (\ref{qi19}), we have
\begin{eqnarray*}
&&Y_{T'-t}^{T'-t,\cdot}(\hat{\omega})={\hat{\theta}}_{T'-t}Y_{0}^{0,\cdot}(\hat{\omega})=Y_{0}^{0,\cdot}({\hat{\theta}}_{T'-t}\hat{\omega}),\\
&&Y_{{T'}^*-t}^{{T'}^*-t,\cdot}({\hat{\omega}}^*)=\hat{\theta}^*_{{T'}^*-t}Y_{0}^{0,\cdot}({\hat{\omega}}^*)=Y_{0}^{0,\cdot}(\hat{\theta}^*_{{T'}^*-t}{\hat{\omega}}^*).
\end{eqnarray*}
So we only need to assert that
${\hat{\theta}}_{T'-t}{\hat{\omega}}=\hat{\theta}^*_{{T'}^*-t}{\hat{\omega}}^*$.
Indeed we have for any $s\geq0$\underline{,}
\begin{eqnarray*}
({\hat{\theta}}_{T'-t}{\hat{\omega}})(s)&=&{\hat{\omega}}(T'-t+s)-{\hat{\omega}}(T'-t)\\
&=&({B}_{T'-(T'-t+s)}-{B}_{T'})-({B}_{T'-(T'-t)}-{B}_{T'})\\
&=&{B}_{t-s}-{B}_{t}.
\end{eqnarray*}
Note that the right hand side of the above formula does not depend
on $T'$, therefore
${\hat{\theta}}_{T'-t}{\hat{\omega}}(s)={\hat{\theta}^*}_{{T'}^*-t}{\hat{\omega}}^*(s)={B}_{t-s}-{B}_{t}$.

On the probability space $(\Omega,\mathscr{F},P)$, we define
${\theta}_{t}=(\hat{\theta}_{t})^{-1}$, $t\geq0$. Actually $\hat{B}$
is a two-sided Brownian motion, so
$(\hat{\theta}_{t})^{-1}=\hat{\theta}_{-t}$ is well defined (see
\cite{ar}). It is easy to see that ${\theta}_{t}$ is a shift w.r.t.
${B}$ satisfying
\begin{description}
\item[$(\textrm{i})$]$P\cdot({\theta}_{t})^{-1}=P$;
\item[$(\textrm{i}\textrm{i})$]${\theta}_{0}=I$;
\item[$(\textrm{i}\textrm{i}\textrm{i})$]${\theta}_{s}\circ{\theta}_{t}={\theta}_{s+t}$;
\item[$(\textrm{iv})$]${\theta}_{t}\circ{B}_s={B}_{s+t}-{B}_{t}$.
\end{description}

Since
$v_t(\cdot)(\omega)=u(T'-t,\cdot)(\hat{\omega})=Y_{T'-t}^{T'-t,\cdot}(\hat{\omega})$
a.s., so
\begin{eqnarray*}
{\theta}_rv_t(\cdot)(\omega)&=&\hat{\theta}_{-r}u(T'-t,\cdot)(\hat{\omega})=\hat{\theta}_{-r}\hat{\theta}_{r}u(T'-t-r,\cdot)(\hat{\omega})\\
&=&u(T'-t-r,\cdot)(\hat{\omega})=v_{t+r}(\cdot)(\omega),
\end{eqnarray*}
for all $r\geq0$ and $T'\geq t+r$ a.s. In particular, let
$Y(\cdot)(\omega)=v_0(\cdot)(\omega)=Y_{T'}^{T',\cdot}(\hat{\omega})$.
Then the above formula implies:
\begin{eqnarray*}
{\theta}_tY(\omega)=Y({\theta}_t\omega)=v_t(\omega)=v(t,v_0(\omega),\omega)=v(t,Y(\omega),\omega)\
{\rm for}\ {\rm all}\ t\geq0\ {\rm a.s.}
\end{eqnarray*}
That is to say
$v_t(\cdot)(\omega)=v_0(\cdot)(\theta_t\omega)=Y(\cdot)({\theta}_t\omega)=Y_{T'-t}^{T'-t,\cdot}(\hat{\omega})$
is a stationary solution of Eq.(\ref{zhao21}) w.r.t. ${\theta}$.
$\hfill\diamond$
\\

\subsection{The solution of infinite horizon BDSDE}
We now consider the following infinite horizon BDSDE with infinite
dimensional noise, which has a more general form than BDSDE
(\ref{qi13}):
\begin{eqnarray}\label{qi30}
{\rm e}^{-Ks}Y_{s}^{t,x}&=&\int_{s}^{\infty}{\rm e}^{-Kr}f(r,X_{r}^{t,x},Y_{r}^{t,x},Z_{r}^{t,x})dr+\int_{s}^{\infty}K{\rm e}^{-Kr}Y_{r}^{t,x}dr\\
&&-\int_{s}^{\infty}{\rm e}^{-Kr}
g(r,X_{r}^{t,x},Y_{r}^{t,x},Z_{r}^{t,x})d^\dagger{\hat{B}}_r-\int_{s}^{\infty}{\rm
e}^{-Kr}\langle Z_{r}^{t,x},dW_r\rangle.\nonumber
\end{eqnarray}
Here
$f:[0,\infty)\times\mathbb{R}^{d}\times\mathbb{R}^1\times\mathbb{R}^{d}{\longrightarrow{\mathbb{R}^1}}$,
$g:[0,\infty)\times\mathbb{R}^{d}\times\mathbb{R}^1\times\mathbb{R}^{d}\longrightarrow
{\mathcal{L}^2_{U_0}(\mathbb{R}^1)}$. Eq.(\ref{qi30}) is equivalent
to
\begin{eqnarray*}\label{qi31}
{\rm e}^{-Ks}Y_{s}^{t,x}&=&\int_{s}^{\infty}{\rm e}^{-Kr}f(r,X_{r}^{t,x},Y_{r}^{t,x},Z_{r}^{t,x})dr+\int_{s}^{\infty}K{\rm e}^{-Kr}Y_{r}^{t,x}dr\nonumber\\
&&-\sum_{j=1}^{\infty}\int_{s}^{\infty}{\rm
e}^{-Kr}g_j(r,X_{r}^{t,x},Y_{r}^{t,x},Z_{r}^{t,x})d^\dagger{\hat{\beta}}_j(r)-\int_{s}^{\infty}{\rm
e}^{-Kr}\langle Z_{r}^{t,x},dW_r\rangle.
\end{eqnarray*}
We assume
\begin{description}
\item[(H.8).] Change ``$\mathscr{B}_{[0,T]}$" to
``$\mathscr{B}_{\mathbb{R}^{+}}$" and ``$r\in[0,T]$" to ``$r\geq0$"
in (H.2).
\item[(H.9).] Change ``$\int_{0}^T$" to ``$\int_{0}^\infty {\rm
e}^{-Kr}$" in (H.3).
\item[(H.10).] Change ``$r\in[0,T]$" to ``$r\geq0$" in (H.4).
\item[(H.11).] Change ``$\mu\in\mathbb{R}^{1}$" to ``$\mu>0$ with $2\mu-K-2C-\sum_{j=1}^{\infty}{C_j}>0$", ``$r\in[0,T]$" to ``$r\geq0$" and ``$\leq\mu\int_{\mathbb{R}^d}U(x){|Y_1(x)-Y_2(x)|}^2\rho^{-1}(x)dx$" to ``$\leq-\mu\int_{\mathbb{R}^d}U(x){|Y_1(x)-Y_2(x)|}^2\rho^{-1}(x)dx$" in
(H.5).
\item[(H.12).] Change ``$r\in[0,T]$" to ``$r\geq0$" in (H.6).
\end{description}
Then we have the existence and uniqueness theorem for the general
form BDSDE (\ref{qi30}):
\begin{thm}\label{zz47} Under Conditions {\rm(H.7)}--{\rm(H.12)},
Eq.(\ref{qi30}) has a unique solution.
\end{thm}
{\em Proof}.
Here we only prove the existence of solution as the uniqueness is
similar to the procedure in the uniqueness proof of Theorem 5.1 in
\cite{zh-zh} although we need the technique as in the uniqueness
proof of Theorem \ref{zz000} to deal with the non-Lipschitz term.
For each $n\in\mathbb{N}$, we define a
sequence of BDSDEs by setting $h=0$ and $T=n$ in Eq.(\ref{qi20}): 
\begin{eqnarray}\label{zz62}
Y_{s}^{t,x,n}&=&\int_{s}^{n}f(r,X_{r}^{t,x},Y_{r}^{t,x,n},Z_{r}^{t,x,n})dr-\int_{s}^{n}g(r,X_{r}^{t,x},Y_{r}^{t,x,n},Z_{r}^{t,x,n})d^\dagger{\hat{B}}_r\nonumber\\
&&-\int_{s}^{n}\langle Z_{r}^{t,x,n},dW_r\rangle,\ \ \ 0\leq s\leq
n.
\end{eqnarray}
It is easy to verify that BDSDE (\ref{zz62}) satisfies conditions of
Theorem \ref{zz48}. Therefore, for each $n$, there exists
$({Y}_\cdot^{t,\cdot,n},{Z}_\cdot^{t,\cdot,n})\in
S^{2,-K}([0,n];L_{\rho}^2({\mathbb{R}^{d}};{\mathbb{R}^{1}}))\bigotimes\\
M^{2,-K}([0,n];L_{\rho}^2({\mathbb{R}^{d}};{\mathbb{R}^{d}}))$ and
$({Y}_s^{t,x,n},{Z}_s^{t,x,n})$ is the unique solution of
Eq.(\ref{zz62}). That is to say, for an arbitrary $\varphi\in
C_c^{0}(\mathbb{R}^d;\mathbb{R}^1)$, $({Y}_s^{t,x,n},{Z}_s^{t,x,n})$
satisfies
\begin{eqnarray}\label{zhang682}
&&\int_{\mathbb{R}^{d}}{\rm e}^{-Ks}Y_s^{t,x,n}\varphi(x)dx\nonumber\\
&=&\int_{s}^{n}\int_{\mathbb{R}^{d}}{\rm e}^{-Kr}f(r,X_{r}^{t,x},Y_{r}^{t,x,n},Z_{r}^{t,x,n})\varphi(x)dxdr\nonumber\\
&&+\int_{s}^{n}\int_{\mathbb{R}^{d}}K{\rm e}^{-Kr}Y_{r}^{t,x,n}\varphi(x)dxdr\nonumber\\
&&-\sum_{j=1}^{\infty}\int_{s}^{n}\int_{\mathbb{R}^{d}}{\rm e}^{-Kr}g_j(r,X_{r}^{t,x},Y_{r}^{t,x,n},Z_{r}^{t,x,n})\varphi(x)dxd^\dagger{\hat{\beta}}_j(r)\nonumber\\
&&-\int_{s}^{n}\langle \int_{\mathbb{R}^{d}}{\rm
e}^{-Kr}Z_{r}^{t,x,n}\varphi(x)dx,dW_r\rangle\ \ \ P-{\rm a.s.}
\end{eqnarray}
Let $(Y_{t}^{n}, Z_{t}^{n})_{t>n}=(0, 0)$, then
$({Y}_\cdot^{t,\cdot,n},{Z}_\cdot^{t,\cdot,n})\in S^{2,-K}\bigcap
M^{2,-K}([0,\infty);L_{\rho}^2({\mathbb{R}^{d}};\\{\mathbb{R}^{1}}))\bigotimes
M^{2,-K}([0,\infty);L_{\rho}^2({\mathbb{R}^{d}};{\mathbb{R}^{d}}))$.
Using a similar argument as in the proof of Theorem 5.1 in
\cite{zh-zh}, we can prove that $({Y}_s^{t,x,n},{Z}_s^{t,x,n})$ is a
Cauchy sequence. Take $({Y}_s^{t,x},{Z}_s^{t,x})$ as the limit of
$({Y}_s^{t,x,n},{Z}_s^{t,x,n})$ in the space $S^{2,-K}\bigcap
M^{2,-K}\\([0,\infty);L_{\rho}^2({\mathbb{R}^{d}};{\mathbb{R}^{1}}))\bigotimes
M^{2,-K}([0,\infty);L_{\rho}^2({\mathbb{R}^{d}};{\mathbb{R}^{d}}))$
and we will show that $({Y}_s^{t,x},{Z}_s^{t,x})$ is the solution of
Eq.(\ref{qi30}). We only need to verify that for arbitrary
$\varphi\in C_c^{0}(\mathbb{R}^d;\mathbb{R}^1)$,
$({Y}_s^{t,x},{Z}_s^{t,x})$ satisfies 
\begin{eqnarray}\label{zz63}
\int_{\mathbb{R}^{d}}{\rm
e}^{-Ks}Y_s^{t,x}\varphi(x)dx&=&\int_{s}^{\infty}\int_{\mathbb{R}^{d}}{\rm e}^{-Kr}f(r,X_{r}^{t,x},Y_r^{t,x},Z_r^{t,x})\varphi(x)dxdr\nonumber\\
&&+\int_{s}^{\infty}\int_{\mathbb{R}^{d}}K{\rm e}^{-Kr}Y_{r}^{t,x}\varphi(x)dxdr\nonumber\\
&&-\sum_{j=1}^{\infty}\int_{s}^{\infty}\int_{\mathbb{R}^{d}}{\rm e}^{-Kr}g_j(r,X_{r}^{t,x},Y_r^{t,x},Z_r^{t,x})\varphi(x)dxd^\dagger{\hat{\beta}}_j(r)\nonumber\\
&&-\int_{s}^{\infty}\langle \int_{\mathbb{R}^{d}}{\rm
e}^{-Kr}Z_r^{t,x}\varphi(x)dx,dW_r\rangle\ \ \ P-{\rm a.s.}
\end{eqnarray}
Noting that $({Y}_s^{t,x,n},{Z}_s^{t,x,n})$ satisfies
Eq.(\ref{zhang682}), we can prove that $({Y}_s^{t,x},{Z}_s^{t,x})$
satisfies Eq.(\ref{zz63}) by verifying that along a subsequence
Eq.(\ref{zhang682}) converges to Eq.(\ref{zz63}) in $L^2(\Omega)$
term by term as $n\longrightarrow\infty$. Here we only show that
along a subsequence
\begin{eqnarray*}
&&E[\ |\int_{s}^{n}\int_{\mathbb{R}^{d}}{\rm
e}^{-Kr}f(r,X_{r}^{t,x},Y_r^{t,x,n},Z_r^{t,x,n})\varphi(x)dxdr\\
&&\ \ \ \ \ \ -\int_{s}^{\infty}\int_{\mathbb{R}^{d}}{\rm
e}^{-Kr}f(r,X_{r}^{t,x},Y_r^{t,x},Z_r^{t,x})\varphi(x)dxdr|^2]\longrightarrow0\
{\rm as}\ n\longrightarrow\infty.
\end{eqnarray*}
For this, note
\begin{eqnarray*}
&&E[\ |\int_{s}^{n}\int_{\mathbb{R}^{d}}{\rm e}^{-Kr}f(r,X_{r}^{t,x},Y_r^{t,x,n},Z_r^{t,x,n})\varphi(x)dxdr\nonumber\\
&&\ \ \ \ \ \ \ -\int_{s}^{\infty}\int_{\mathbb{R}^{d}}{\rm e}^{-Kr}f(r,X_{r}^{t,x},Y_r^{t,x},Z_r^{t,x})\varphi(x)dxdr|^2]\nonumber\\
&\leq&2E[\ |\int_{s}^{n}\int_{\mathbb{R}^{d}}{\rm e}^{-Kr}\big(f(r,X_{r}^{t,x},Y_r^{t,x,n},Z_r^{t,x,n})\nonumber\\
&&\ \ \ \ \ \ \ \ \ \ \ \ \ \ \ \ \ \ \ \ \ \ \ \ \ -f(r,X_{r}^{t,x},Y_r^{t,x},Z_r^{t,x})\big)\varphi(x)dxdr|^2]\nonumber\\
&&+2E[\ |\int_{n}^{\infty}\int_{\mathbb{R}^{d}}{\rm e}^{-Kr}f(r,X_{r}^{t,x},Y_r^{t,x},Z_r^{t,x})\varphi(x)dxdr|^2]\nonumber\\
&\leq&C_pE[\int_{s}^{\infty}\int_{\mathbb{R}^{d}}{\rm e}^{-Kr}|f(r,X_{r}^{t,x},Y_r^{t,x,n},Z_r^{t,x,n})\nonumber\\
&&\ \ \ \ \ \ \ \ \ \ \ \ \ \ \ \ \ \ \ \ \ \ \ \ \ -f(r,X_{r}^{t,x},Y_r^{t,x,n},Z_r^{t,x})|^2\rho^{-1}(x)dxdr]\nonumber\\
&&+C_pE[\int_{s}^{n}\int_{\mathbb{R}^{d}}{\rm e}^{-Kr}|f(r,X_{r}^{t,x},Y_r^{t,x,n},Z_r^{t,x})\nonumber\\
&&\ \ \ \ \ \ \ \ \ \ \ \ \ \ \ \ \ \ \ \ \ \ \ \ \ \ -f(r,X_{r}^{t,x},Y_r^{t,x},Z_r^{t,x})|^2\rho^{-1}(x)dxdr]\nonumber\\
&&+C_pE[\int_{n}^{\infty}\int_{\mathbb{R}^{d}}{\rm e}^{-Kr}|f(r,X_{r}^{t,x},Y_r^{t,x},Z_r^{t,x})|^2\rho^{-1}(x)dxdr]\nonumber\\
&\leq&C_pE[\int_{s}^{\infty}\int_{\mathbb{R}^{d}}{\rm e}^{-Kr}|Z_r^{t,x,n}-Z_r^{t,x}|^2\rho^{-1}(x)dxdr]\nonumber\\
&&+C_pE[\int_{n}^{\infty}\int_{\mathbb{R}^{d}}{\rm e}^{-Kr}(1+|Y_r^{t,x}|^2+|Z_r^{t,x}|^2)\rho^{-1}(x)dxdr]\nonumber\\
&&+C_pE[\int_{s}^{\infty}\int_{\mathbb{R}^{d}}{\rm
e}^{-Kr}|f(r,X_{r}^{t,x},Y_r^{t,x,n},Z_r^{t,x})\nonumber\\
&&\ \ \ \ \ \ \ \ \ \ \ \ \ \ \ \ \ \ \ \ \ \ \ \ \ \ \
-f(r,X_{r}^{t,x},Y_r^{t,x},Z_r^{t,x})|^2\rho^{-1}(x)dxdr].
\end{eqnarray*}
Similar to (\ref{zz17}), we only need to prove that along a
subsequence
\begin{eqnarray}\label{zhangb48}
&&E[\int_{s}^{\infty}\int_{\mathbb{R}^{d}}{\rm
e}^{-Kr}|f(r,X_{r}^{t,x},Y_r^{t,x,n},Z_r^{t,x})\\
&&\ \ \ \ \ \ \ \ \ \ \ \ \ \ \ \ \ \ \ \ \ \
-f(r,X_{r}^{t,x},Y_r^{t,x},Z_r^{t,x})|^2\rho^{-1}(x)dxdr]\longrightarrow0\
{\rm as}\ n\longrightarrow\infty.\nonumber
\end{eqnarray}
Since $\{{Y}_s^{t,x,n}\}_{n=1}^\infty$ is a Cauchy sequence in the
space
$M^{2,-K}([0,\infty);L_{\rho}^2({\mathbb{R}^{d}};{\mathbb{R}^{1}}))$
with the limit ${Y}_s^{t,x}$, as $n\rightarrow0$, we have
\begin{eqnarray}\label{zz19}
E[\int_{0}^{\infty}\int_{\mathbb{R}^{d}}{\rm
e}^{-Kr}|{Y}_{r}^{t,x,n}-Y_r^{t,x}|^2\rho^{-1}(x)dxdr]\longrightarrow0.
\end{eqnarray}
Then from (\ref{zz19}) we can find a subsequence of
$\{Y_r^{t,x,n}\}_{n=1}^\infty$ still denoted by
$\{Y_r^{t,x,n}\}_{n=1}^\infty$ s.t. $Y_r^{t,x,n}\longrightarrow
Y_r^{t,x}$ for a.e. $r\geq0$, $x\in\mathbb{R}^{d}$, a.s. $\omega$
and $E[\int_{0}^{\infty}\int_{\mathbb{R}^{d}}{\rm
e}^{-Kr}\sup_n|Y_r^{t,x,n}|^2\rho^{-1}(x)dxdr]<\infty$. Therefore,
for this subsequence $\{Y_r^{t,x,n}\}_{n=1}^\infty$, by Condition
(H.10), we have
\begin{eqnarray*}
&&E[\int_{0}^{\infty}\int_{\mathbb{R}^{d}}{\rm e}^{-Kr}\sup_n|f(r,X_{r}^{t,x},Y_r^{t,x,n},Z_r^{t,x})\nonumber\\
&&\ \ \ \ \ \ \ \ \ \ \ \ \ \ \ \ \ \ \ \ \ \ \ \ \ \ \ -f(r,X_{r}^{t,x},Y_r^{t,x},Z_r^{t,x})|^2\rho^{-1}(x)dxdr]\\
&\leq&
C_pE[\int_{0}^{\infty}\int_{\mathbb{R}^{d}}{\rm e}^{-Kr}(1+\sup_n|Y_r^{t,x,n}|^2+|Y_r^{t,x}|^2+|Z_r^{t,x}|^2)\rho^{-1}(x)dxdr]\\
&<&\infty.
\end{eqnarray*}
Then (\ref{zhangb48}) follows from applying Lebesgue's dominated
convergence theorem and Condition (H.12). 
That is to say ${(Y_s^{t,x}, Z_s^{t,x})}_{s\geq0}$ satisfies
Eq.(\ref{zz63}). The proof of Theorem \ref{zz47} is completed.
$\hfill\diamond$\\

By a similar method as in the proof of the existence part in case
(i) in Theorem 5.1 in \cite{zh-zh}, we have the following
estimation:
\begin{prop}\label{qi072}
Let $({Y}_s^{t,x,n},{Z}_s^{t,x,n})$ be the solution of
Eq.(\ref{zz62}), then under the conditions of Theorem \ref{zz47},
\begin{eqnarray*}
&&\sup_nE[\sup_{s\geq0}\int_{\mathbb{R}^d}{\rm
e}^{-Ks}|Y_s^{t,x,n}(x)|^2\rho^{-1}(x)dx]\nonumber\\
&&+\sup_nE[\int_{0}^{\infty}\int_{\mathbb{R}^{d}}{\rm
e}^{-Kr}|Y_r^{t,x,n}(x)|^2\rho^{-1}(x)dxdr]\nonumber\\
&&+\sup_nE[\int_{0}^{\infty}\int_{\mathbb{R}^{d}}{\rm
e}^{-Kr}|Z_r^{t,x,n}(x)|^2\rho^{-1}(x)dxdr]<\infty.
\end{eqnarray*}
\end{prop}

\subsection{Proofs of Theorem \ref{qi043} and Theorem \ref{qi044}}
All the proofs until now in this paper have shown us how to deal
with the non-Lipschitz term. Indeed the proofs of Theorem
\ref{qi043} and Theorem \ref{qi044} are rather similar to the proofs
in Section 6 in \cite{zh-zh} even under the non-Lipschitz
conditions. So we only intend to give the proof briefly.

{\em Proof of Theorem \ref{qi043}}. Since the conditions here are
stronger than those in Theorem \ref{zz47}, so there exists a unique
solution $(Y_{s}^{t,x},Z_{s}^{t,x})$ to Eq.(\ref{qi13}). We only
need to prove $E[\sup_{s\geq0}\int_{\mathbb{R}^{d}}{\rm
e}^{{-{pK}}s}{{|{Y}_s^{t,x}|}^p}\rho^{-1}(x)dx]<\infty$. Let
$\varphi_{N,p}(x)=x^{p\over2}I_{\{0\leq
x<N\}}+N^{{p-2}\over2}({p\over2}x-{{p-2}\over2}N)I_{\{x\geq N\}}$.
We apply the generalized It$\hat {\rm o}$'s formula to ${\rm
e}^{{-{pK}}r}\varphi_{N,p}\big(\psi_M(Y_{r}^{t,x})\big)$ for a.e.
$x\in{\mathbb{R}^{d}}$ to have
\begin{eqnarray}\label{zhang690}
&&{\rm e}^{{-{pK}}s}\varphi_{N,p}\big(\psi_M(Y_{s}^{t,x})\big)-{pK}\int_{s}^{T}{\rm e}^{{-{pK}}r}\varphi_{N,p}\big(\psi_M(Y_{r}^{t,x})\big)dr\nonumber\\
&&+{1\over2}\int_{s}^{T}{\rm e}^{{-{pK}}r}\varphi^{''}_{N,p}\big(\psi_M(Y_{r}^{t,x})\big)|\psi_M^{'}(Y_{r}^{t,x})|^2|Z_{r}^{t,x}|^2dr\nonumber\\
&&+\int_{s}^{T}{\rm e}^{{-{pK}}r}\varphi^{'}_{N,p}\big(\psi_M(Y_{r}^{t,x})\big)I_{\{-M\leq{Y}_{r}^{t,x}<M\}}|{Z}_{r}^{t,x}|^2dr\nonumber\\
&\leq&{\rm e}^{{-{pK}}T}\varphi_{N,p}\big(\psi_M(Y_{T}^{t,x})\big)\nonumber\\
&&+\int_{s}^{T}{\rm e}^{{-{pK}}r}\varphi^{'}_{N,p}\big(\psi_M(Y_{r}^{t,x})\big)\psi_M^{'}(Y_{r}^{t,x}){f}(X_{r}^{t,x},Y_{r}^{t,x},Z_{r}^{t,x})dr\nonumber\\
&&+\int_{s}^{T}{\rm e}^{{-{pK}}r}\varphi^{'}_{N,p}\big(\psi_M(Y_{r}^{t,x})\big)I_{\{-M\leq{Y}_{r}^{t,x}<M\}}\sum_{j=1}^{\infty}|{g}_j(X_{r}^{t,x},Y_{r}^{t,x},Z_{r}^{t,x})|^2dr\nonumber\\
&&+{1\over2}\int_{s}^{T}{\rm e}^{{-{pK}}r}\varphi^{''}_{N,p}\big(\psi_M(Y_{r}^{t,x})\big)|\psi_M^{'}(Y_{r}^{t,x})|^2\sum_{j=1}^{\infty}|{g}_j(X_{r}^{t,x},Y_{r}^{t,x},Z_{r}^{t,x})|^2dr\nonumber\\
&&-\sum_{j=1}^{\infty}\int_{s}^{T}{\rm e}^{{-{pK}}r}\varphi^{'}_{N,p}\big(\psi_M(Y_{r}^{t,x})\big)\psi_M^{'}(Y_{r}^{t,x}){g}_j(X_{r}^{t,x},{Y}_{r}^{t,x},{Z}_{r}^{t,x})d^\dagger{\hat{\beta}}_j(r)\nonumber\\
&&-\int_{s}^{T}\langle{\rm
e}^{{-{pK}}r}\varphi^{'}_{N,p}\big(\psi_M(Y_{r}^{t,x})\big)\psi_M^{'}(Y_{r}^{t,x}){Z}_{r}^{t,x},dW_r\rangle.
\end{eqnarray}
From the above estimation, using $\lim_{T\rightarrow\infty}{\rm
e}^{{-{pK}}T}\varphi_{N,p}\big(\psi_M(Y_{T}^{t,x})\big)=0$ and
taking the limit as $M\to \infty$ first, then the limit as $N\to
\infty$, by the monotone convergence theorem, we have
\begin{eqnarray}\label{zhang691}
&&E[\int_{s}^{\infty}\int_{\mathbb{R}^{d}}{\rm e}^{{-{pK}}r}{|{Y}_r^{t,x}|}^p\rho^{-1}(x)dxdr]\nonumber\\
&&+E[\int_{s}^{\infty}\int_{\mathbb{R}^{d}}{\rm e}^{{-{pK}}r}{{|{Y}_r^{t,x}|}^{p-2}}|{Z}_{r}^{t,x}|^2\rho^{-1}(x)dxdr]\nonumber\\
&\leq&C_p+C_p\int_{\mathbb{R}^{d}}\sum_{j=1}^{\infty}|{g}_j(x,0,0)|^p\rho^{-1}(x)dx]<\infty.
\end{eqnarray}
Also by the B-D-G inequality, the Cauchy-Schwartz inequality and the
Young inequality, we can obtain another estimation from
(\ref{zhang690}):
\begin{eqnarray*}\label{zhang692}
&&E[\sup_{s\geq0}\int_{\mathbb{R}^{d}}{\rm e}^{{-{pK}}s}{{|{Y}_s^{t,x}|}^p}\rho^{-1}(x)dx]\nonumber\\
&\leq&C_p\int_{\mathbb{R}^{d}}|{f}(x,0,0)|^p\rho^{-1}(x)dx+C_p\int_{\mathbb{R}^{d}}\sum_{j=1}^{\infty}|{g}_j(x,0,0)|^p\rho^{-1}(x)dx\\
&&+C_pE[\int_{0}^{\infty}\int_{\mathbb{R}^{d}}{\rm
e}^{{-{pK}}r}{{|{Y}_r^{t,x}|}^{p-2}}|{Z}_{r}^{t,x}|^2\rho^{-1}(x)dxdr]\nonumber\\
&&+C_pE[\int_{0}^{\infty}\int_{\mathbb{R}^{d}}{\rm
e}^{{-{pK}}r}|Y_{r}^{t,x}|^p\rho^{-1}(x)dxdr].\nonumber
\end{eqnarray*}
So by (\ref{zhang691}), Theorem \ref{qi043} is proved.
$\hfill\diamond$
\\

Now we turn to the proof of Theorem \ref{qi044}.

{\em Proof of Theorem \ref{qi044}}. First note that we also can
prove Lemma 6.2 in \cite{zh-zh} under the conditions in this
theorem, so we have
\begin{eqnarray*}
&&E([\sup_{s\geq0}\int_{\mathbb{R}^{d}}{\rm e}^{-{2K}s}|Y_{s}^{t',x}-Y_{s}^{t,x}|^2\rho^{-1}(x)dx])^{p\over2}\nonumber\\
&\leq&C_pE[\sup_{s\geq0}\int_{\mathbb{R}^{d}}{\rm e}^{-{pK}r}|Y_{s}^{t',x}-Y_{s}^{t,x}|^p\rho^{-1}(x)dx]\big(\int_{\mathbb{R}^{d}}\rho^{-1}(x)dx\big)^{{p-2}\over2}\nonumber\\
&\leq&C_p|t'-t|^{p\over2}.
\end{eqnarray*}
This is because we actually did not use the Lipschitz condition of
$f$ w.r.t. $y$ and the monotone condition is enough. Noting $p>2$,
by the Kolmogorov continuity theorem (see \cite{ku2}), we have
$t\longrightarrow Y_{s}^{t,x}$ is a.s. continuous for $t\in[0,T]$
under the norm $(\sup_{s\geq0}\int_{\mathbb{R}^{d}}{\rm
e}^{-{2K}s}|\cdot|^2\rho^{-1}(x)dx)^{1\over2}$. Without losing any
generality, assume that $t'\geq t$. Then we can see that
\begin{eqnarray*}
&&\lim_{t'\rightarrow t}(\int_{\mathbb{R}^{d}}{\rm
e}^{-{2K}t'}|Y_{t'}^{t',x}-Y_{t'}^{t,x}|^2\rho^{-1}(x)dx)^{1\over2}\\
&\leq&\lim_{t'\rightarrow t}(\sup_{s\geq0}\int_{\mathbb{R}^{d}}{\rm
e}^{-{2K}s}|Y_{s}^{t',x}-Y_{s}^{t,x}|^2\rho^{-1}(x)dx)^{1\over2}=0\
\ \ {\rm a.s.}
\end{eqnarray*}
Notice $t'\in[0,T]$, so
\begin{eqnarray}\label{zhang706}
\lim_{t'\rightarrow
t}(\int_{\mathbb{R}^{d}}|Y_{t'}^{t',x}-Y_{t'}^{t,x}|^2\rho^{-1}(x)dx)^{1\over2}=0\
\ \ \ {\rm a.s.}
\end{eqnarray}
Since $Y_{\cdot}^{t,\cdot}\in
S^{2,-K}([0,\infty);L_{\rho}^2({\mathbb{R}^{d}};{\mathbb{R}^{1}}))$,
$Y_{t'}^{t,\cdot}$ is continuous w.r.t. $t'$ in
$L_{\rho}^2(\mathbb{R}^d;\mathbb{R}^1)$. That is to say for each
$t$,
\begin{eqnarray}\label{zhang705}
\lim_{t'\rightarrow
t}(\int_{\mathbb{R}^{d}}|{Y}_{t'}^{t,x}-{Y}_{t}^{t,x}|^2\rho^{-1}(x)dx)^{1\over2}=0\
\ \ \ {\rm a.s.}
\end{eqnarray}
Now by (\ref{zhang706}) and (\ref{zhang705})
\begin{eqnarray*}
&&\lim_{t'\rightarrow t}(\int_{\mathbb{R}^{d}}|{Y}_{t'}^{t',x}-{Y}_{t}^{t,x}|^2\rho^{-1}(x)dx)^{1\over2}\\
&\leq&\lim_{t'\rightarrow t}(\int_{\mathbb{R}^{d}}|{Y}_{t'}^{t',x}-{Y}_{t'}^{t,x}|^2\rho^{-1}(x)dx)^{1\over2}+\lim_{t'\rightarrow t}(\int_{\mathbb{R}^{d}}|{Y}_{t'}^{t,x}-{Y}_{t}^{t,x}|^2\rho^{-1}(x)dx)^{1\over2}\\
&=&0\ \ \ \ {\rm a.s.}
\end{eqnarray*}
For arbitrary $T>0$, $0\leq t\leq T$, define
$u(t,\cdot)=Y_t^{t,\cdot}$, then $u(t,\cdot)$ is a.s. continuous
w.r.t. $t$ in $L_{\rho}^2(\mathbb{R}^d;\mathbb{R}^1)$. Since
$Y_{\cdot}^{t,\cdot}\in
S^{2,-K}([0,\infty);L_{\rho}^2({\mathbb{R}^{d}};{\mathbb{R}^{1}}))$,
$Y_T^{T,x}$ is
$\mathscr{F}_{T,\infty}^{\hat{B}}\otimes\mathscr{B}_{\mathbb{R}^{d}}$
measurable and
$E[\int_{\mathbb{R}^{d}}|Y_T^{T,x}|^2\rho^{-1}(x)dx]<\infty$. It
follows that Condition (H.1) is satisfied. Moreover, Conditions
{\rm(A.1)}--{\rm(A.6)} are stronger than Conditions (H.2)--(H.7), so
by Theorem \ref{zz006}, $u(t,x)$ is a weak solution of
Eq.(\ref{zhang685}). Theorem \ref{qi044} is proved. $\hfill\diamond$
\\
\\
{\bf Acknowledgements}. It is our great pleasure to thank S. Peng
and K. Lu for useful conversations. QZ would like to acknowledge the
financial support of the National Basic Research Program of China
(973 Program) with Grant No. 2007CB814904.


\begin{thebibliography}{[99]}
\footnotesize

\bibitem{ar} L. Arnold, {\em Random dynamical systems}. Springer-Verlag Berlin Heidelberg (1998).

\bibitem{ar-sc} L. Arnold, M. Scheutzow, {\em Perfect cocycles through stochastic differential equations}. Probab. Theory Relat. Fields, Vol.101 (1995), 65-88.

\bibitem{ba-ma} V. Bally, A. Matoussi, {\em Weak solutions for SPDEs and backward doubly stochastic differential equations}. Journal of Theoretical Probability, Vol.14 (2001), 125-164.

\bibitem{ba-le} G. Barles, E. Lesigne, {\em SDE, BSDE and PDE}. In: Backward stochastic differential equations. Pitman Res. Notes Math., Ser.364, Longman, Harlow, (1997), 47-80.

\bibitem{br-de-hu-pr-st} Ph. Briand, B. Delyon, Y. Hu, E. Pardoux, L. Stoica, {\em $L^p$ solutions of backward stochastic differential equations}. Stochastic Process. Appl., Vol.108 (2003), 109-129.

\bibitem{br-hu} Ph. Briand, Y. Hu, {\em BSDE with quadratic growth and unbounded terminal value}. Probab. Theory Relat. Fields, Vol.136 (2006), 604-618.



\bibitem{kloeden} T. Caraballo, P.E. Kloeden, B. Schmalfuss, {\em Exponentially stable stationary solutions for stochastic evolution equations and their perturbation}. Appl. Math. Optim., Vol.50 (2004), 183--207.

\bibitem{pr-za1} G. Da Prato, J. Zabczyk, {\em Stochastic equations in infinite dimensions}. Cambridge University Press (1992).


\bibitem{du-lu-sc1} J. Duan, K. Lu, B. Schmalfuss, {\em Invariant manifolds for stochastic partial differential equations}. Ann. Probab., Vol.31 (2003), 2109-2135.





\bibitem{el-tr-zh} K. D. Elworthy, A. Truman, H. Z. Zhao, {\em Generalized It$\hat {o}$ formulae and space-time Lebesgue-Stieltjes integrals of local times}. S${\rm\grave{e}}$minaire de Probabilit${\rm\grave{e}}$s, Vol.40 (2007), 117-136.

\bibitem{ha} R. Z. Has$'$minskii, {\em Stochastic stability of differential equations}. Alphen aan den Rijn (The Netherlands): Sijthoff and Noordhoff (1980).

\bibitem{ko} M. Kobylanski, {\em Backward stochastic differential equations and partial differential equations with quadratic growth}. Ann. Probab., Vol.28, No.2 (2000), 558-602.

\bibitem{krylov} N.V. Krylov, {\em An analytic approach to SPDEs}, in: Stochastic partial equations:
six perspectives, edited by R.A. Carmona and B. Rozovskii,
Mathematical Surveys and Monographs, Vol. 64, American Mathematical
Society (1999), 183-242.

\bibitem{ku2} H. Kunita, {\em Stochastic flows and stochastic differential equations}. Cambridge University Press (1990).

\bibitem{ku1} H. Kunita, {\em Stochastic flow acting on Schwartz distributions}. J. Theor. Prob., 7(2), (1994), 247-278.

\bibitem{le-sa} J. P. Lepeltier, J. San Martin, {\em Backward stochastic differential equations with continuous coefficient}. Statistics and Probability Letters, Vol.32 (1997), 425-430.



\bibitem{mo-zh-zh} S.-E. A. Mohammed, T. Zhang, H. Z. Zhao, {\em The stable manifold theorem for semilinear stochastic evolution equations and stochastic partial differential equations}. Memoirs of the American Mathematical Society, Vol.196 (2008), No. 917, 105pp .

\bibitem{pa} E. Pardoux, {\em Backward stochastic differential equations and viscosity solutions of systems of semilinear parabolic and elliptic PDEs of second order}. Stochastic Analysis and Related Topics: The Geilo Workshop, 1996, L. Decreusefond, J. Gjerde, B. Oksendal, A.S. Ust¨¹unel eds., Birkh$\ddot{a}$user, (1998), 79-127.


\bibitem{pa-pe1} E. Pardoux, S. Peng, {\em Adapted solution of a backward stochastic differential equation}. Syst. Control Lett., Vol.14 (1990), 55-61.


\bibitem{pa-pe3} E. Pardoux, S. Peng, {\em Backward doubly stochastic differential equations and systems of quasilinear SPDEs}. Probab. Theory Relat. Fields, Vol.98 (1994), 209-227.

\bibitem{pe} S. Peng, {\em Probabilistic interpretation for systems of quasilinear parabolic partial differential equations}. Stochastics, Vol.37 (1991), 61-74.


\bibitem{shi} Y. F. Shi, Y. L. Gu, K. Liu, {\em Comparison theorems of backward doubly stochastic differential equations and applications}. Stochastic Analysis and Applications, Vol.23 (2005), 97-110.

\bibitem{si1} Ya. Sinai, {\em Two results concerning asymptotic behaviour of solutions of Burgers equation with force}. J. Statist. Phys., Vol.64 (1991), 1-12.

\bibitem{si2} Ya. Sinai, {\em Burgers system driven by a periodic stochastic flows}. In: It$\hat {\rm o}$'s stochastic calculus and probability theory, Springer, Tokyo, (1996), 347-353.

\bibitem{zh-zh} Q. Zhang, H. Z. Zhao, {\em Stationary solutions of SPDEs and infinite horizon BDSDEs}. Journal of Functional Analysis, Vol.252 (2007), 171-219.
\end{thebibliography}
\end{document}